\documentclass[a4paper,11pt,openright,twoside]{article}

\usepackage{color}
\usepackage{latexsym}
\usepackage{amssymb,amsbsy,amsmath,amsfonts,amssymb,amscd}
\usepackage{epsfig, graphicx}
\usepackage{hyperref}
\usepackage{enumerate}
\usepackage{mathrsfs}
\usepackage{color}
\usepackage{delarray}

\usepackage{float}
\usepackage{tikz}
\usepackage[all]{xy}
\usepackage{wrapfig}
\usepackage{enumitem}
\usepackage{multirow, array} % para las tablas
\usepackage{setspace}
\usepackage{appendix}
\usepackage{fancyhdr}
\newcommand{\commentout}[1]{}
\newcommand{\R}{\mathbb{R}}
\newcommand{\N}{\mathbb{N}}

\newcommand {\e}  {\mathrm{e}}

\def\tu{\widetilde u}

\def\tn{\widetilde  n }

\newcommand{\norme}[1]{\left\lVert#1\right\rVert}

\newcommand{\pscal}[2]{\left\langle ~#1~ , ~#2~\right\rangle}

\newcommand{\beq}{\begin{equation}}
\newcommand{\eeq}{\end{equation}}
\newcommand{\bea} {\begin{array}{rl}}
\newcommand{\eea} {\end{array}}
\newcommand{\bepa}{\left\{ \begin{array}{l}}
\newcommand{\eepa} {\end{array}\right.}
\newtheorem{theorem}{Theorem}[section]
\newtheorem{lemma}[theorem]{Lemma}
\newtheorem{definition}[theorem]{Definition}
\newtheorem{remark}[theorem]{Remark}

\newtheorem{hyp}[theorem]{Hypothesis}

\numberwithin{equation}{section}
\newcommand{\qed}{{ \hfill
                      {\unskip\kern 6pt\penalty 500 \raise -2pt\hbox{\vrule\vbox to 6pt{\hrule width 6pt
                      \vfill\hrule}\vrule} \par}   }}

\title{{\Large \bf A hyperbolic-elliptic-parabolic PDE model \\
	describing chemotactic {\it E.~coli} colonies}\\
	~\\
	 {\large Cet article est dédié à la mémoire du Professeur Ezzedine Zahrouni.}}
\author{Danielle Hilhorst\thanks{CNRS, Laboratoire de Math\'{e}matiques d'Orsay (LMO), Universit\'{e} Paris-Saclay, Orsay, France.
	\texttt{danielle.hilhorst@math.u-psud.fr}}
\and Pierre Roux\thanks{Laboratoire de Math\'{e}matiques d'Orsay (LMO), Universit\'{e} Paris-Saclay, Orsay, France. 	\texttt{pierre.roux@math.u-psud.fr}}
}

\date{\today}

\begin{document}
\maketitle           
\begin{abstract}
  We study a modified version of an initial-boundary value problem describing the formation of colony patterns of bacteria \textit{Escherichia Coli}. The original system of three parabolic equations was studied numerically and analytically and gave insights into the underlying mechanisms of chemotaxis. We focus here on the parabolic-elliptic-parabolic approximation and the hyperbolic-elliptic-parabolic limiting system which describes the case of pure chemotactic movement without random diffusion. We first construct local-in-time solutions for the parabolic-elliptic-parabolic system. Then we prove uniform \textit{a priori} estimates and we use them along with a compactness argument in order to construct local-in-time solutions for the hyperbolic-elliptic-parabolic limiting system. Finally, we prove that some initial conditions give rise to solutions which blow-up in finite time in the $L^\infty$ norm in all space dimensions. This last result violet is true even in space dimension 1, which is not the case for the full parabolic or parabolic-elliptic Keller-Segel systems.
\end{abstract}

\vspace{5cm}

{\bf Key-words}: Blow-up, Chemotaxis, Keller-Segel, Hyperbolic equation.
\\
\\
{\bf AMS Class. No}: 35B44, 35A01, 35L04 %35K60, 82C31, 92B20

\pagenumbering{arabic}

%%%%%%%%%%%%%%%%%%%%%%%

%%%%%%%%%%%%%%%%%%%%%%%%%%%%%%%%%%%%%%%%%%%%%%%%%%%%%5

\newpage
\section{Introduction}

\setlength{\belowdisplayskip}{4pt} \setlength{\belowdisplayshortskip}{4pt}
\setlength{\abovedisplayskip}{4pt} \setlength{\abovedisplayshortskip}{4pt}

Since the 70' and the first article of Keller and Segel \cite{KS2} taxis-diffusion and aggregation equations were widely studied and proved to model a wide range of phenomena. This first model was of the form
\begin{equation}\label{KS}
\left\{\begin{array}{rcl}
\partial_t u&=&\nabla \cdot ( D_2(u,v) \nabla u)  - \nabla \cdot (D_1(u,v)\nabla v ),\\
\partial_t v&=& D \Delta v + f(u)u - k(v)v
\end{array} \right.
\end{equation}
and aimed at understanding the aggregation of amoebae (namely the species \textit{Dictyostelium discoideum}). Here $u$ denote the density of the amoebae and $v$ the density of a chemotactic substance which is produced by the amoebae in case of food scarceness. Many variants were introduced and mathematically studied. Some of these models are subject to blow-up in finite time (see the monograph \cite{B1} and references therein). We refer to \cite{H1}, \cite{H2}, \cite{B1}, \cite{HP} and \cite{BBTW} for reviews of the ever-growing bibliography on this subject.
Some of these variants are growth models: the first equation is modified to allow for birth and death of organisms. Several growth terms were studied in early works (for example \cite{BHLM1},\cite{BHLM2},\cite{HHS},\cite{FMT},\cite{OTYM},\cite{FHMW},\cite{MT},\cite{OTYM}. The Keller-Segel system with logistic growth $ru - \mu u^2$ was intensively studied, with the recent contributions by Tello and Winkler \cite{TW}, Winkler \cite{W1}, Lankeit \cite{L1}, \cite{L2}, \cite{L3}, Kang and Stevens \cite{KS1}.

Among other generalisations with growth terms, Aotani, Mimura and Mollee \cite{AMM} developed a Keller-Segel type model that reproduced qualitatively the experimental results obtained by Budrene and Berg \cite{BB1}, \cite{BB2}, where bacteria produce stable geometrical patterns under the effect of chemotaxis. This new model aimed to investigate the hypothesis of \cite{BB1} that cell proliferation, excretion of chemoattractant, chemotactic motility and substrate consumption can generate the observed patterns without any specialised morphogenetic program. The numerical simulations of \cite{AMM} gave insights in this direction. This model was investigated in dimension 1 in \cite{HMT}. The model is given by
\begin{equation}\label{Mim}
\left\{\begin{array}{rcl}
\partial_t u&=&\Delta u - \nabla \cdot (u\nabla \chi (c) ) + g(u)nu -b(n)u,\\
\partial_t c&=&d_c \Delta c + \alpha u - \beta c,\\
\partial_t n&=&d_n \Delta n - \gamma g(u)n u,\\
\partial_t w&=&b(n)u,\\
u(x,0)&=&u_0(x), \ c(x,0)=c_0(x), \ n(x,0)=n_0(x), \ w(x,0) = w_0(x),
\end{array} \right.
\end{equation}
where we denote by $u(x,t)$ the density of active bacteria, $w(x,t)$ the density of inactive bacteria, $c(x,t)$ the concentration of chemoattractant and $n(x,t)$ the density of nutrients. The constants $d_u, d_n, \alpha, \beta, \gamma$ are assumed to be positive and $g,b$ and $\chi$ are regular enough functions which we will discuss below. The system is studied for $t\in[0,+\infty)$ in a open, bounded domain $\Omega\subset \R^d$ with smooth boundary and Neumann boundary conditions:
\begin{equation} \dfrac{\partial u}{\partial \nu}=\dfrac{\partial c}{\partial \nu}=\dfrac{\partial n}{\partial \nu}=0 \qquad \qquad  \mathrm{on} \ \partial\Omega. \end{equation}

Celinski, Karch, Hilhorst, Mimura and Roux \cite{CHKMR} prove the existence and uniqueness of local-in-time solutions and construct global-in-time solutions in the one dimensional case. Under smallness assumptions on the initial data, they also establish a global-in-time existence result in higher dimension. In the case that the solution $(u,c,n,w)$ exists globally in time, it converges exponentially fast in $L^\infty$ norm, as $t$ tends to $+\infty$, towards an unknown stationary pattern of the form $(0,0,n_\infty,w_\infty(\cdot))$ where $n_\infty$ is a non-negative constant and $w_\infty$ a non-negative function. For a simplified parabolic-elliptic-parabolic (see problem \eqref{MimuraEps} below) version of the system, the authors prove a blow-up result for initial data with large enough $L^1$ norm and sufficiently concentrated around a given point.\\

Following the articles of Winkler \cite{W1}, Lankeit \cite{L3} and Kang and Stevens \cite{KS1} in the logistic growth case, we study here a hyperbolic-elliptic-parabolic version of the model proposed in \cite{AMM} where the diffusion term is removed from the first equation:
\begin{equation}\label{Intermediaire}
\left\{\begin{array}{rcll}
\partial_t u&=& - \nabla \cdot (u\nabla c) + g(u)nu - b(n)u,\qquad &\mathrm{in} \ \Omega\times\R^+\\
0&=&\Delta c + \alpha u - \beta c,&\mathrm{in} \ \Omega\times\R^+\\
\partial_t n&=&\Delta n - \gamma g(u)n u,&\mathrm{in} \ \Omega\times\R^+\\
\dfrac{\partial c}{\partial \nu} &=& \dfrac{\partial n}{\partial \nu}\ =\ 0 & \mathrm{on} \ \partial\Omega\times\R^+.\\
u(x,0)&=&u_0(x), \qquad n(x,0)=n_0(x)&\mathrm{in} \ \Omega.
\end{array} \right.
\end{equation}
in an open, smooth, and bounded domain $\Omega\subset\R^d$.
%\begin{equation}\label{Bound} \end{equation}
We also define
\begin{equation} \label{Bound}
	w(x,t)=\int_0^t b\big(n(x,s)\big)u(x,s) ds, \quad \mathrm{in} \ \Omega\times\R^+
\end{equation}
as the density of inactive bacteria for this system.\\

We assume that
\begin{itemize}
	\item[$\bullet$] $g$ is a $\mathscr{C}^1$ positive increasing function such that $g(0)=0$, $\lim_{x\to+\infty}g(x)=G_0$ and $g'\in L^\infty(\R^+)$;
	\item[$\bullet$] $b$ is a $\mathscr{C}^1$ positive decreasing function such that $b(0)=B_0>0$ and $b'\in L^\infty(\R^+)$;
	\item[$\bullet$] $\alpha,\beta,\gamma$ are positive constants.
\end{itemize}
 and that the initial conditions satisfy the following hypothesis:
\begin{hyp}\label{H0}
		$u_0,n_0\in W^{2,p}(\Omega)$ for some $p>d$, $u_0,n_0\geqslant 0$ and $u_0,n_0$ satisfy homogeneous Neumann boundary conditions.
\end{hyp}
 
\begin{remark}
 	Since the boundary $\partial \Omega$ is a characteristic, no boundary condition is needed for $u$.
\end{remark}
 
This system models the case where random motion of bacteria is negligible compared to the effect of chemotaxis. It depicts well the case of a high density of chemoattractant in the environment when bacteria already colonized the medium and start gathering and inactivating themselves. Note that hyperbolic Keller-Segel systems were the subject of active research in the past few years (see \cite{PD}, \cite{FGM} and bibliography in \cite{BBTW}).

In Section 2, we first construct local-in-time solutions for the parabolic-elliptic-parabolic approximate system \eqref{MimuraEps}. Then, under the assumption that the domain $\Omega$ is convex, we prove that the solution of an ODE is an upper-solution of the $W^{1,p}$ norm of the component $u_\epsilon$ of solutions $(u_\varepsilon,c_\varepsilon,n_\varepsilon)$ of \eqref{MimuraEps} on a time interval $[0,T]$, where $T$ is independent from the diffusion coefficient $\varepsilon$ in the first equation of \eqref{MimuraEps}. Hence, we deduce from the comparison principle a $W^{1,p}$ bound on the solutions in $[0,T]$ which is uniform with respect to $\varepsilon$. Other estimates follow by Sobolev injection and computations. In Section 3, we assume again that the space domain $\Omega$ is convex and we construct local-in-time solutions for system \eqref{Intermediaire} \textit{via} a compactness method involving the previous \textit{a priori} estimates on the approximate system. The time of existence of these solutions is at least $T$, where $T$ only depends on the norm of the initial condition. In Section 4, we prove a blow-up result for solutions of \eqref{Intermediaire} for all space dimensions, including dimension 1. In the appendix A, we recall some classical results which we need all along : technical lemmas, elliptic estimates and heat semi-group estimates.

%%%%%%%%%%%%%%%%%%%%%%%%%%%%%%%%%%%%%%%%%%%%%%

\newpage

\section{The parabolic-elliptic-parabolic approximate system}
\label{sec: local}

Let $\varepsilon>0$ be arbitrary. In order to construct local-in-time solutions to system \eqref{Intermediaire}, we consider the approximate system
\begin{equation}\label{MimuraEps}
\left\{\begin{array}{rcll}
\partial_t u_\varepsilon&=&\varepsilon\Delta u_\varepsilon - \nabla \cdot (u_\varepsilon\nabla c_\varepsilon ) + g(u_\varepsilon)n_\varepsilon u_\varepsilon - b(n_\varepsilon)u_\varepsilon,\qquad &\mathrm{in} \ \Omega\times\R^+\\
0&=&\Delta c_\varepsilon + \alpha u_\varepsilon - \beta c_\varepsilon, &\mathrm{in} \ \Omega\times\R^+\\
\partial_t n_\varepsilon&=&\Delta n_\varepsilon - \gamma g(u_\varepsilon)n_\varepsilon u_\varepsilon, &\mathrm{in} \ \Omega\times\R^+\\
\dfrac{\partial u_\varepsilon}{\partial \nu} &=& \dfrac{\partial c_\varepsilon}{\partial \nu}\ =\ \dfrac{\partial n_\varepsilon}{\partial \nu}\ =\ 0 &  \mathrm{on} \ \partial\Omega\times\R^+\\
u_\varepsilon(x,0)&=&u_0^\varepsilon(x), \qquad n_\varepsilon(x,0)=n_0^\varepsilon(x) &\mathrm{in} \ \Omega,
\end{array} \right.
\end{equation}
where $u_0^\varepsilon,n_0^\varepsilon$ satisfy the hypothesis:
\begin{hyp}\label{H0eps}
	$u_0^\varepsilon,n_0^\varepsilon\in C^{2,\eta}(\bar \Omega)$ for some $\eta\in(0,1)$, $u_0^\varepsilon,n_0^\varepsilon\geqslant 0$ and $u_0^\varepsilon,n_0^\varepsilon$ satisfy homogeneous Neumann boundary conditions. Moreover, there exists $p>d$ such that for all $\varepsilon>0$,
	\[ \norme{u_0^\varepsilon}_{W^{p,\infty}} \leqslant \norme{u_0}_{W^{p,\infty}} ,\qquad \norme{n_0^\varepsilon}_{W^{1,\infty}} \leqslant \norme{n_0}_{W^{1,\infty}} \qquad \mathrm{and} \qquad   \norme{n_0^\varepsilon}_\infty \leqslant \norme{n_0}_\infty . \] 
	When $\varepsilon$ tends to $0$, $u_0^\varepsilon$ and $n_0^\varepsilon$ converge strongly in $W^{2,p}(\Omega)$ towards $u_0$ and $n_0$ respectively, where $u_0$ and $n_0$ satisfy Hypothesis \ref{H0}.
\end{hyp}

\subsection{Local-in-time existence}

Consider the heat semigroup $\{e^{t\Delta}\}_{t\geq 0}$ in the smooth bounded domain $\Omega$ with Neumann boundary conditions. 
We can write the following Duhamel's formulae for problem \eqref{MimuraEps}:
\begin{align}\label{duhE1}
\begin{split}
u_\varepsilon(t)=&\ {\rm e}^{\varepsilon\Delta t} u_0^\varepsilon +\int_0^t{\rm e}^{\varepsilon\Delta(t-s)} \nabla\cdot \Big( u_\varepsilon(s)\nabla c_\varepsilon(s) \Big) d s
+\int_0^t{\rm e}^{\varepsilon\Delta(t-s)}u_\varepsilon(s)(g(u_\varepsilon)n_\varepsilon-b(n_\varepsilon))(s) d s,
\end{split}
\\
%\begin{align}
n_\varepsilon(t)=&\ {\rm e}^{\Delta t} n_0^\varepsilon -\gamma\int_0^t{\rm e}^{\Delta(t-s)}g(u_\varepsilon(s))n_\varepsilon(s)u_\varepsilon(s) d s.\label{duhE2}
\end{align}
where $c_\varepsilon(t)$ is the unique solution of
\begin{equation}\label{duhE3}
\begin{array}{l}
-\Delta c_\varepsilon(t) +\beta c_\varepsilon(t) = \alpha u_\varepsilon(t)  \quad \text{for}\quad x\in\Omega\quad\text{and}\quad t>0,\\
\dfrac{\partial c_\varepsilon}{\partial \nu}=0 \quad  \text{for}\quad x\in\partial\Omega\quad\text{and}\quad t>0.
\end{array}
\end{equation}

\begin{definition}
	Let $\Omega\subset \R^d$ an open domain with smooth boundary.
%	\begin{itemize}
%		\item Let $u:\bar \Omega\to\R$, we say that $u\in C^{0}(\bar \Omega)$ if $u$ is bounded and continuous and we denote
%		\[  \norme{u}_\infty = \sup_{x\in\Omega} |u(x)|.  \]
%		\item Let $\eta\in(0,1)$. Let $u\in C^{0}(\bar \Omega)$, we say that $u\in C^{0,\eta}(\bar \Omega)$, if
%		\[ \norme{u}_{C^{0,\eta}} = \norme{u}_\infty + \sup_{x,y\in\Omega, x\neq y} \left\{  \dfrac{|u(x)-u(y)|}{|x-y|^\eta}   \right\}<+\infty.  \]
%		\item Let $k\in\N$ and $\eta\in(0,1)$. Let $u\in C^{k}(\bar \Omega)$, we say that $u\in C^{k,\eta}(\bar \Omega)$ if
%		\[  \norme{u}_{C^{k,\alpha}} = \sum_{|\alpha|\leqslant k} \norme{ \partial^\alpha u  }_\infty + \sum_{|\alpha|=k} \sup_{x,y\in\Omega, x\neq y} \left\{  \dfrac{|\partial^\alpha u(x)-\partial^\alpha u(y)|}{|x-y|^\eta}   \right\}<+\infty. \]
		We define
			\[\mathcal{V}_1 = \{  u\in C^{\eta,\frac\eta2}(\bar\Omega\times[0,T]) \ | \ \nabla u \in C^{\eta,\frac\eta2}(\bar\Omega\times[0,T]) \}\]
		and
		\[\mathcal{V}_2 = \{  u\in \mathcal{V}_1 \ | \ \nabla^2 u \in C^{\eta,\frac\eta2}(\bar\Omega\times[0,T]) \}.\]
%		\item Let $\eta\in(0,1]$, $k,m\in\N$ and $T\in[0,+\infty)$. We say that $u:\bar \Omega\times[0,T)\to\R$, $u\in C^{k+\eta_1,m+\eta_2}(\bar \Omega\times[0,T))$ if for all $\tilde T\in(0,T)$, $u\in C^{k+\eta_1,m+\eta_2}(\bar \Omega\times[0,\tilde T])$.  We say that $u:\bar \Omega\times[0,T)\to\R$, $u\in C^{k+\eta_1,m+\eta_2}(\bar \Omega\times(0,T))$ if for all $T_1,T_2\in(0,T)$ such that $T_1<T_2$, $u\in C^{k+\eta_1,m+\eta_2}(\bar \Omega\times[T_1,T_2])$
%	\end{itemize}
\end{definition}

%\cap C^{1+\eta,\frac\eta2}(\bar\Omega\times[0,T_{max}))
%[Celinski, Karch, Hilhorst, Mimura, Roux]
\begin{theorem}\label{thm:localE}
	Let $u_0^\varepsilon,n_0^\varepsilon$ satisfy Hypothesis \ref{H0eps}. There exists a unique maximal time $T_{max,\varepsilon}\in(0,+\infty]$ such that problem 
	\eqref{MimuraEps} has a unique non-negative solution 
	$$
	(u_\varepsilon,c_\varepsilon,n_\varepsilon) \in C^{2+\eta,1+\frac\eta2}(\bar\Omega\times[0,T_{\mathrm{max},\varepsilon})),
	$$
	which satisfies for all $t\in[0,T_{\mathrm{max},\varepsilon})$, $ \norme{n_\varepsilon(\cdot,t)} \leqslant \norme{ n_0^\varepsilon }_\infty $.
	Moreover, if $T_{\mathrm{max},\varepsilon}<+\infty$, then
	\begin{equation}\label{Boom} \limsup_{t\to T_{\mathrm{max},\varepsilon}}\norme{u_\varepsilon(t)}_{\infty}=+\infty. \end{equation}
\end{theorem}

\noindent{\bf Proof.}
	First, by Lemma \ref{ElliptReg1} if $u\in C^{\eta,\frac\eta2}(\bar\Omega\times[0,T])$, then the solution $c_\varepsilon$ of \eqref{duhE3} satisfies $c_\varepsilon\in \mathcal{V}_2$ and for all $t>0$
	\begin{equation} \label{gradEst} \norme{ \nabla c_\varepsilon(t) }_\infty \leqslant  K_\infty\norme{ u_\varepsilon(t) }_\infty, \end{equation}
	for some constant $K_\infty>0$.\\
	
	Then, we fix $R>0$ such that
	\[ \max(\norme{  u_0^\varepsilon  }_\infty,\norme{  n_0^\varepsilon  }_\infty) \leqslant \dfrac{R}{2},  \]
	and we define the space
	\begin{equation*}
	V_{T,R}= \big\{ u,n \in\, \mathcal{V}_1, \ 
	\sup_{0\leqslant t\leqslant T} \norme{u(t)}_\infty\leqslant R,\quad 
	\sup_{0\leqslant t\leqslant T} \norme{n(t)}_\infty\leqslant R
	\big\}.
	\end{equation*}
	Let $\mathcal{F}= \mathcal{F}(u_\varepsilon,n_\varepsilon): V_{T,R}\to V_{T,R}$ be the mapping defined by the right-hand side of equations \eqref{duhE1}, \eqref{duhE2} and \eqref{duhE3}.
	
	Let $(u_\varepsilon,n_\varepsilon)\in V_{T,R}$. we define	$	(\tu,\tn)	$	by the right-hand sides of \eqref{duhE1} and \eqref{duhE2} applied to $u_\varepsilon$ and $n_\varepsilon$, where $c_\varepsilon$ satisfies \eqref{duhE3}. These holds:
	\begin{multline*}
	\norme{  \tu(t) }_\infty \ \leq\  \norme{ e^{\varepsilon\Delta t} u_0^\varepsilon  }_\infty  +  \norme{\  \int_0^t{\rm e}^{\varepsilon \Delta(t-s)}\nabla\cdot\Big(u_\varepsilon(s)\nabla c_\varepsilon(s)   \Big)  d s }_\infty \\
	+  \norme{ \int_0^t{\rm e}^{\varepsilon \Delta(t-s)}\Big(u_\varepsilon(s)(g(u_\varepsilon)n_\varepsilon-b(n_\varepsilon))(s)\Big) d s
	\ }_\infty
	\end{multline*}
	and
	\begin{equation*}
	\norme{  \tn(t) }_\infty \ \leq\   \norme{ e^{t\Delta} n_0^\varepsilon }_\infty +   \gamma\norme{ \int_0^t{\rm e}^{\Delta(t-s)}\Big(g(u_\varepsilon(s))n_\varepsilon(s)u_\varepsilon(s)\Big) d s  }_\infty.
	\end{equation*}
	Moreover, in view of \eqref{gradEst}
	\begin{equation*}
	\norme{u_\varepsilon\nabla  c_\varepsilon }_\infty	\leqslant  \norme{ u_\varepsilon(t) }_\infty \norme{ \nabla c_\varepsilon }_\infty \leqslant K_\infty R^2,
	\end{equation*}
	\begin{equation*}
	\norme{\big(g(u_\varepsilon)n_\varepsilon - b(n_\varepsilon)\big)u_\varepsilon}_\infty\leqslant \big(G_0\norme{n_\varepsilon}_\infty + B_0\big) \norme{u_\varepsilon}_\infty \leqslant (G_0R + B_0\big)R.
	\end{equation*}
	We combine these inequalities with the estimates of the heat semigroup
	\begin{equation*}
	\begin{split}
	\norme{e^{t\Delta} v}_\infty \leqslant \norme{v}_\infty \quad \text{and}\quad 
	\norme{ e^{t\Delta} \nabla\cdot w}_\infty \leqslant C\left(1+ t^{-\frac12}\right)\norme{w}_\infty
	\end{split}
	\end{equation*}
	(see Lemma~\ref{lem:LpLq} below) to deduce that
	\begin{align}
	\sup_{0\leqslant t\leqslant T} \norme{  \tu(t) }_\infty \leqslant &\, \dfrac{R}{2} +  (T+2\sqrt{T})CK_\infty  R^2 + T(G_0R + B_0\big)R, \label{I1}  \\
	\sup_{0\leqslant t\leqslant T} \norme{  \tn(t) }_\infty \leqslant &\,  \dfrac{R}{2} + T\gamma  G_0R^2 \label{I2}
	\end{align}
	Furthermore, by Lemma \ref{ElliptReg1}, if $u\in C^{\eta,\frac\eta2}(\bar\Omega\times[0,T])$, then the solution $c$ of \eqref{duhE3} satisfies $c\in \mathcal{V}_2$. Hence, $\nabla  c  \in \mathcal{V}_1$ and thus $u\nabla c  \in \mathcal{V}_1$. Then, by the definition and properties of $\{e^{\Delta t}\}_{t\geq 0}$, the semigroup of linear operators
	generated by the Laplacian with the Neumann boundary condition, and by the properties of the integral, we have
	\[{\rm e}^{\varepsilon\Delta(t-s)}\nabla \cdot \Big(u(s)\nabla c(s) \Big)\in C^{2+\eta,1+\frac\eta2}(\bar\Omega\times[0,T]),\]
	and thus
	\[  t\mapsto \int_0^t{\rm e}^{\varepsilon\Delta(t-s)}\nabla\cdot\Big(u(s)\nabla c(s)   \Big)  d s  \in C^{2+\eta,1+\frac\eta2}(\bar\Omega\times[0,T]). \]
	Applying these arguments to each term in the right-hand sides of the Duhamel's formulae \eqref{duhE1}, \eqref{duhE2} and \eqref{duhE3}, we conclude, thanks to Hypothesis \ref{H0eps}, that
	\begin{equation}
	\tu,\tn\in C^{2+\eta,1+\frac\eta2}(\bar\Omega\times[0,T]) \subset \mathcal{V}_1.
	\end{equation}
	We choose $T=T(R)$ depending on $R$ small enough so that the right-hand sides of the inequalities \eqref{I1} and \eqref{I2} are bounded from above by $R$. Then, we have $\mathcal{F} ( V_{T,R} ) \subset V_{T,R}$.
	Since the embedding $ C^{2+\eta,1+\frac\eta2}(\bar\Omega\times[0,T]) \hookrightarrow \mathcal{V}_1 $ is compact, $\mathcal{F} ( V_{T,R} )$ is a relatively compact subset of $V_{T,R}$. Hence, the continuous mapping $\mathcal{F}$ has a fixed point by the Schauder fixed point theorem.\\
	
	The non-negativity of $u_\varepsilon$, $c_\varepsilon$, $n_\varepsilon$ is a direct consequence of the maximum principle for parabolic and elliptic equations. The maximum principle applied to the equation for $n_\varepsilon$ also yields $\norme{n_\varepsilon(\cdot,t)}_\infty\leqslant \norme{ n_0^\varepsilon }_\infty$.
	
	We now prove uniqueness : assume there exist two solutions $(u_1,c_1,n_1)$ and $(u_2,c_2,n_2)$ on $[0,T]$. Let $U,C,N$ be defined by $U=u_1-u_2$, $C=c_1-c_2$ and $N=n_1-n_2$.\\
	We choose $ \mathcal R$ large enough such that
	\[  \sup_{t\in[0,T]} \norme{ u_i(t)}_{W^{1,\infty}} \leqslant  \mathcal R  , \ \sup_{t\in[0,T]} \norme{ c_i(t)}_{W^{1,\infty}} \leqslant  \mathcal R  , \ \sup_{t\in[0,T]} \norme{ n_i(t)}_\infty \leqslant  \mathcal R   \]
	for  $i=1,2$.
	Since, using the equation for $c_i$, $i=1,2$,
	\[ \Delta c_i = \beta c_i-\alpha u_i, \]
	we have, for $i=1,2$,
	\begin{equation}\label{Pari}  \dfrac{\partial u_i}{\partial t} - \varepsilon\Delta u_i =  - \nabla  u_i \cdot \nabla c_i - u_i(\beta c_i-\alpha u_i)   + (g(u_i)n_i - b(n_i))u_i.  \end{equation}
	Hence,
	\begin{multline}\label{Par}  \dfrac{\partial U}{\partial t} - \varepsilon\Delta U =  - \nabla   U\cdot \nabla c_1 - \nabla  u_2 \cdot \nabla C  + \Big(g(u_1)n_1 - b(n_1)   - (\beta c_1-\alpha u_1)  \Big)u_1\\ - \Big(g(u_2)n_2-b(n_2)  - (\beta c_2-\alpha u_2)\Big)u_2.  \end{multline}
	We have
	\begin{equation*}
	|b(n_1)u_1-b(n_2)u_2|\leqslant   
	|U|\sup_{|n|\leqslant  \mathcal R  }|b(n)|
	+|N|\sup_{s\in[0,T]}\norme{u_2(s)}_\infty \sup_{|n|\leqslant  \mathcal R  }|b'(n)|,
	\end{equation*}
	\begin{equation*}
	|g(u_1)u_1n_1-g(u_2)u_2n_2| \leqslant 
	|N |\sup_{|u|\leqslant  \mathcal R} |g(u)u| +
	|U|\sup_{s\in[0,T]}\norme{n_2(s)}_\infty 
	\sup_{|u|\leqslant  \mathcal R} |g'(u)u+g(u)| ,
	\end{equation*}
	\begin{equation*} | u_2(\beta c_2-\alpha u_2) - u_1(\beta c_1-\alpha u_1) | \leqslant \beta | u_2c_2 - u_1c_1 | + \alpha | u_2^2 - u_1^2 |   \\
	  \leqslant \beta \mathcal R (  |C| + |U| ) +  2\alpha \mathcal R |U|  \end{equation*} 
	Before multiplying \eqref{Par} by $U$, we give the further estimates: 
	\[ |  \nabla   U\cdot \nabla c_1  U  | \leqslant  \frac\varepsilon 2| \nabla   U |^2    +   \frac1{2\varepsilon}(  \sup_{s\in[0,T]} \norme{\nabla c_1(s)}_\infty)^2U^2 \]
	\begin{equation*} | \nabla  u_2 \cdot( \nabla c_1 -\nabla c_2 ) U |   \leqslant  \frac12 (\sup_{s\in[0,T]} \norme{ \nabla u_2(s)}_\infty)^2U^2   +   \frac12 | \nabla C |^2.   \end{equation*}
	Hence, multiplying equation \eqref{Par} by $U$ and integrating, we obtain
	\begin{equation} \label{ineqU} \dfrac{d}{dt} \int_\Omega U^2  + \varepsilon \int_ \Omega |\nabla U |^2    \leqslant  \kappa_1(\varepsilon) \Big(  \int_\Omega U^2 + \int_\Omega C^2 + \int_\Omega | \nabla C |^2 +\int_\Omega N^2   \Big)  , \end{equation}
	for some large enough constant $\kappa_1=k_1(\varepsilon)>0$.\\
	Then, we have
	\[  -\Delta C + \beta C = \alpha U. \]
	We multiply by $C$ and we integrate in order to obtain,
	\[   \int_\Omega | \nabla C |^2  +   \beta \int_\Omega C^2 = \alpha \int_\Omega C U \leqslant \alpha\Big(   \dfrac{\beta}{2\alpha} \int_\Omega C^2   + \dfrac{\alpha}{2\beta}    \int_\Omega U^2    \Big)  ,  \]
	and as a consequence
	\begin{equation} \label{ineqC} \int_\Omega | \nabla C |^2 + \frac{\beta}{2} \int_\Omega C^2 \leqslant \frac{\alpha^2}{2\beta} \int_\Omega U^2.  \end{equation}
	Last, we have
	\[  \dfrac{\partial N}{\partial t} = \Delta N - \gamma( g(u_1)n_1u_1 - g(u_2)n_2u_2 ). \]
	We multiply by $N$ and we integrate. With the same reasoning as above we find that, for some large enough constant $\kappa_2>0$,
	\begin{equation}\label{IneqN}
	\dfrac{d}{dt}\int_\Omega N^2 + 2\int_{\Omega} |\nabla N|^2 \leqslant \kappa_2 \Big(  \int_\Omega N^2 + \int_\Omega U^2   \Big)
	\end{equation}
	
	We deduce from \eqref{ineqU}, \eqref{ineqC} and \eqref{IneqN} that, for some large enough constant $\lambda>0$,
	\begin{equation} \dfrac{d}{dt} \int_\Omega (U^2 + N^2 ) \leqslant \lambda \int_\Omega (U^2 + N^2), \end{equation}
	and -- since $U(0)^2 + N(0)^2=0$ -- applying Gronwall's lemma gives $U(t)=0$ and $N(t)=0$ for all $t\in[0,T]$. As a consequence, $C(t)=0$ for all $t\in[0,T]$. This completes the proof of the uniqueness of the solution.\\
	
	Finally, we argue by contradiction in order to prove \eqref{Boom}. Assume that the maximal time of existence $T_{\mathrm{max},\varepsilon}$ for the unique solution is finite. Assume that
	\[ \max\left(  \sup_{t\in[0,T_{\mathrm{max},\varepsilon})} \norme{ u_\varepsilon(t) }_\infty, \norme{ n_0^\varepsilon }_\infty \right) < \frac{\tilde R_\varepsilon}2 < +\infty.  \]
	Then, in the construction of local-in-time solutions done above, we can choose a constant $\tilde R_\varepsilon$ that works for every $(u_\varepsilon(t),n_\varepsilon(t))$ taken as an initial condition. As stated above, the local time of existence $T(\tilde R_\varepsilon)$ from any initial datum $(u_\varepsilon(t),n_\varepsilon(t))$, $t<T_{\mathrm{max},\varepsilon}$ only depends on $\tilde R_\varepsilon$ and other constants. Hence, we can choose $t$ close enough to $T_{\mathrm{max},\varepsilon}$ so as to achieve $t+T(\tilde R_\varepsilon) > T_{\mathrm{max},\varepsilon}$. It's a contradiction with $T_{\mathrm{max},\varepsilon}$ being the maximal time of existence. Therefore,
	\[ \sup_{t\in[0,T_{\mathrm{max},\varepsilon})} \norme{ u_\varepsilon(t) }_\infty = +\infty.  \]
\qed

\subsection{Uniform \textit{a priori} estimates}

We derive \textit{a priori} estimates in the spirit of \cite{KS1}.

\begin{lemma}\label{lemme2}
	Assume $\Omega$ is bounded and convex, with a smooth boundary. Let $u_0^\varepsilon,n_0^\varepsilon$ satisfy hypothesis \ref{H0eps}. There exists a time $T>0$ and there exist positive constants $M_1=M_1(T)$, $M_2=M_2(T)$ depending on $T$ such that for all $\varepsilon>0$, the solution $(u_\varepsilon,c_\varepsilon,n_\varepsilon)$ of \eqref{MimuraEps} constructed in Theorem \ref{thm:localE} exists in $\Omega\times[0,T]$ and we have \begin{equation} \label{L2B1} \sup_{t\in[0,T]}\norme{u_\varepsilon(t)}_{\infty} < M_1(T),\end{equation} and 
	\begin{equation} \label{L2B2} \sup_{t\in[0,T]}\norme{\nabla u_\varepsilon(t)}_{p} < M_2(T).\end{equation}
\end{lemma}

\noindent{\bf Proof.}
Let $p\in(d,+\infty)$; note that we have the Sobolev embedding $W^{1,p}(\Omega)\hookrightarrow L^\infty(\Omega) $, so that, for all $f\in W^{1,p}(\Omega)$,
\begin{equation} \label{Embedd}
	\norme{f}_\infty \leqslant C_{\infty,p} \norme{f}_{W^{1,p}},
\end{equation}
with $C_{\infty,p}>0$ depending only on $\Omega$ and $p$ (cf. \cite{B2}, Corollary 9.14).

First, we multiply the first equation of \eqref{MimuraEps} by $u_\varepsilon^{p-1}$ and we compute
\begin{equation*} \dfrac{1}{p} \dfrac{d}{dt} \int_{\Omega} u_\varepsilon^p + \dfrac{4\varepsilon (p-1)}{p^2}\int_{\Omega} \left| \nabla u_\varepsilon^{\frac{p}{2}} \right|^2 \leqslant \alpha\dfrac{p-1}{p}\int_{\Omega}u_\varepsilon^{p+1} -\beta\frac{p-1}{p}\int_{\Omega}c_\varepsilon u_\varepsilon^{p}   + \int_{\Omega} (g(u_\varepsilon)n_\varepsilon-b(n_\varepsilon))u_\varepsilon^p. \end{equation*}
Using \eqref{Embedd}, we have
\[ \int_{\Omega}  u_\varepsilon^{p+1} \leqslant   \norme{u_\varepsilon}_\infty \int_{\Omega}u_\varepsilon^{p} \leqslant   C_{\infty,p} \norme{u_\varepsilon}_{W^{1,p}} \norme{u_\varepsilon}_p^p.  \]
Using elliptic regularity (Lemma \ref{ElliptReg1}), we have
\[ \int_{\Omega}c_\varepsilon u_\varepsilon^{p} \leqslant K_\infty   \norme{u_\varepsilon}_\infty\int_{\Omega} u_\varepsilon^{p} \leqslant K_\infty C_{\infty,p}   \norme{u_\varepsilon}_{W^{1,p}} \norme{u_\varepsilon}^p_p .\]
We deduce from the previous results that there exists a constant $C_1>0$, depending on the parameters of the system, $\Omega$ and $p$, independent from $\varepsilon$, such that
\begin{equation} \label{Brique1}
\dfrac{d}{dt} \norme{ u_\varepsilon }_p^p \leqslant C_1(1 + \norme{u_\varepsilon}_{W^{1,p}} ) \norme{u_\varepsilon}_p^p.
\end{equation}

Second, we denote $w_k=\dfrac{\partial u_\varepsilon}{\partial x_k}$ and noting that
\[ \begin{array}{rcl} \nabla \cdot (u_\varepsilon\nabla c_\varepsilon) &=& \nabla u_\varepsilon \cdot \nabla c_\varepsilon + u_\varepsilon\Delta c_\varepsilon\\
&=& \nabla u_\varepsilon \cdot \nabla c_\varepsilon +  u_\varepsilon(\beta c_\varepsilon-\alpha u_\varepsilon) ,
\end{array} \]
we compute
\begin{multline}\label{Horrible} 
\dfrac{\partial w_k}{\partial t} - \varepsilon \Delta w_k\,=\, \displaystyle   -\,\nabla w_k \cdot \nabla c_\varepsilon -\, \nabla u_\varepsilon \cdot \nabla \dfrac{\partial c_\varepsilon}{\partial x_k}  +\left( g(u_\varepsilon)\dfrac{\partial n_\varepsilon}{\partial x_k} -\beta\dfrac{\partial c_\varepsilon}{\partial x_k}  -b'(n_\varepsilon)\dfrac{\partial n_\varepsilon}{\partial x_k} \right)u_\varepsilon\\
 +\Big(2\alpha u_\varepsilon - \beta c_\varepsilon + g'(u_\varepsilon)n_\varepsilon u_\varepsilon +g(u_\varepsilon)n_\varepsilon - b(n_\varepsilon)\Big)w_k.
 \end{multline}
We denote 
\[ \mathcal{A}(k) = \dfrac{1}{p}\dfrac{d \norme{w_k}_p^p}{d t} - \varepsilon\int_\Omega |w_k|^{p-2}w_k\Delta w_k .  \]
Multiplying \eqref{Horrible} by $|w_k|^{p-2}w_k$ and integrating, we obtain
\begin{equation} \begin{array}{rcl}
\displaystyle \mathcal{A}(k) &=& \displaystyle \frac{1}{p}\int_\Omega  \Delta c_\varepsilon|w_k|^p -\int_\Omega   \nabla u_\varepsilon\cdot \nabla \dfrac{\partial c_\varepsilon}{\partial x_k}|w_k|^{p-2}w_k \\
&&\displaystyle \quad \quad \quad +\int_\Omega \left( g(u_\varepsilon)\dfrac{\partial n_\varepsilon}{\partial x_k}  -\beta \dfrac{\partial c_\varepsilon}{\partial x_k}  -b'(n_\varepsilon)\dfrac{\partial n_\varepsilon}{\partial x_k} \right)u_\varepsilon|w_k|^{p-2}w_k\\
&&\displaystyle \quad \quad \quad \quad +\int_\Omega \Big(2\alpha u_\varepsilon - \beta c_\varepsilon+ g'(u_\varepsilon) n_\varepsilon u_\varepsilon+g(u_\varepsilon)n_\varepsilon - b(n_\varepsilon)\Big)|w_k|^p.
\end{array} \end{equation}
Now, let us derive a bound for $\mathcal{A}(k)$ in terms of $L^\infty$ and $W^{1,p}$ norms of $u$.
In view of hypotheses on $g,b$ and using $\norme{n(\cdot,t)}\leqslant \norme{ n_0^\varepsilon }\leqslant \norme{n_0}_\infty$, \eqref{est2} and $|\Delta c_\varepsilon|\leqslant (\alpha + \beta K_\infty) \norme{u_\varepsilon}_\infty$, we deduce that there exists a large enough constant $C_2>0$ such that 
\begin{equation} \begin{array}{rcl}
\displaystyle \mathcal{A}(k)&\leqslant& \displaystyle C_2\norme{u_\varepsilon}_\infty\int_\Omega  |w_k|^p +C_2 \int_\Omega  |\nabla u_\varepsilon| | \nabla ^2c_\varepsilon||w_k|^{p-1}\\
&&\displaystyle \quad \quad \quad +C_2\big(\norme{u_\varepsilon}_\infty+ \norme{u_\varepsilon}_\infty^2\big)\int_\Omega \Big( G_0\left|\dfrac{\partial n_\varepsilon}{\partial x_k}\right|  + 1   +\norme{b'}_\infty\left|\dfrac{\partial n_\varepsilon}{\partial x_k}\right| \Big)|w_k|^{p-1}\\
&&\displaystyle \quad \quad \quad \quad +C_2\Big(  \norme{u_\varepsilon}_\infty + \norme{g'}_\infty \norme{n_0}_\infty \norme{u_\varepsilon}_\infty + G_0 \norme{n_0}_\infty + B_0\Big)\int_\Omega  |w_k|^p.
\end{array} \end{equation}
The function $n$ satisfies the Duhamel's formulae
\[ n_\varepsilon(x,t)=\e^{\Delta t} n_0^\varepsilon (x,t)-\gamma\int_0^t \e^{\Delta (t-s)}g(u_\varepsilon(s,t))n_\varepsilon(s,t)u_\varepsilon(s,t) d s  \]
and
\[  \nabla n_\varepsilon(x,t) = \nabla \e^{\Delta t} n_0^\varepsilon -\gamma\int_0^t \nabla \Big( \e^{\Delta (t-s)}g(u_\varepsilon(s,t))n_\varepsilon(s,t)u_\varepsilon(s,t) \Big) d s. \]
By Lemma \ref{lem:LpLq} and noting that for all $t>0$,
\[  \norme{\nabla (\e^{\Delta t}  n_0^\varepsilon  )}_\infty \leqslant \norme{\nabla  n_0^\varepsilon  }_\infty \leqslant \norme{ n_0^\varepsilon }_{W^{1,\infty}} ,\]
we have
\[ \begin{array}{rcl} \displaystyle \norme{\nabla n_\varepsilon}_\infty &\leqslant& \displaystyle \norme{ n_0^\varepsilon }_{W^{1,\infty}} + \gamma C \int_0^t \left(1+ (t-s)^{-\frac12}\right)\e^{-\lambda_1(t-s)}\norme{(g(u_\varepsilon)n_\varepsilon u_\varepsilon)(\cdot,s)}_\infty d s \\
&\leqslant& \displaystyle \norme{ n_0^\varepsilon }_{W^{1,\infty}} + \gamma CG_0\norme{ n_0^\varepsilon }_\infty \int_0^t \left( 1 + (t-s)^{-\frac12} \right)\e^{-\lambda_1(t-s)}\norme{u_\varepsilon}_\infty d s\\
&\leqslant& \displaystyle \norme{ n_0^\varepsilon }_{W^{1,\infty}} + \gamma CG_0\norme{ n_0^\varepsilon }_\infty \norme{u_\varepsilon}_\infty \left( t + \int_0^t \frac{1}{\sqrt{t-s}} d s \right) \\
&\leqslant& \displaystyle \norme{n_0}_{W^{1,\infty}} + \gamma CG_0\norme{n_0}_\infty \norme{u_\varepsilon}_\infty\left( t + 2\sqrt{t}\right). \\
\end{array}  \]
Next, we define
\begin{equation}\label{function_f}
	\mathcal F(t) = 1 + t + 2\sqrt{t}
\end{equation}
Hence, there exists a constant $C_3>0$ depending on $\norme{n_0}_{W^{1,\infty}},\Omega,G_0$ and $\gamma$ such that
\begin{equation}\label{NightWork}\norme{\dfrac{\partial n_\varepsilon}{\partial x_k}}_\infty\leqslant C_3\mathcal F(t)(1+ \norme{u_\varepsilon}_\infty ).\end{equation}
Therefore, using Lemma \ref{ElliptReg1}, there exists a constant $C_4 > 0$ such that
\begin{equation} 
\displaystyle \mathcal{A}(k)\leqslant \displaystyle C_4\mathcal F(t)\Big(1+ \norme{u_\varepsilon}_\infty + \norme{u_\varepsilon}_\infty^2  \Big) \Big( \int_\Omega  |\nabla u_\varepsilon| | \nabla ^2c_\varepsilon||w_k|^{p-1}   + \int_\Omega  |w_k|^p + \int_\Omega |w_k|^{p-1} \Big).
\end{equation}
By Lemmas \ref{ElliptReg1}, \ref{BMO1} and \ref{BMO2}, there exist constants $\tilde C_5,C_5>0$ such that
\begin{equation}\begin{array}{rcl}
 \displaystyle \norme{\nabla^2c_\varepsilon}_\infty &\leqslant&\displaystyle \tilde C_5 \left(1+\norme{\nabla^2c_\varepsilon}_{BMO} \left( 1+ \log^+(\norme{\nabla^2c_\varepsilon}_{W^{1,p}}) \right) \right) \\
 &\leqslant&\displaystyle  C_5 \left(1+\norme{u_\varepsilon}_\infty \left( 1+ \log^+(\norme{ u_\varepsilon}_{W^{1,p}}) \right) \right).
\end{array} \end{equation}
Collecting the previous estimates, we deduce that there exists a constant $C_6>0$ depending on the parameters of the system, on $p$, $\norme{n_0}_{W^{1,\infty}}$, $\Omega$ and independent from $\varepsilon$, such that
\begin{equation}\label{presk1} \sum_{k=1}^d \mathcal{A}(k) \leqslant C_6\mathcal F(t) \Big(1 + \norme{u_\varepsilon}_{\infty}^3(1+\log^+(\norme{ u_\varepsilon}_{W^{1,p}}))\Big)\norme{\nabla u_\varepsilon }_p^p. \end{equation}
Moreover, using integration by parts, we have
\begin{equation*} \begin{array}{rcl} \displaystyle \int_\Omega |w_k|^{p-2}w_k\Delta w_k &=&\displaystyle  - (p-1)\int_{\Omega} |\nabla w_k|^2 |w_k|^{p-2}  + \int_{\partial\Omega} \dfrac{\partial w}{\partial \nu} w |w|^{p-2} d\sigma \\
&=&\displaystyle  - (p-1)\int_{\Omega} \left(  |w_k|^{\frac p 2-1} \nabla w_k   \right) \cdot \left(  |w_k|^{\frac p 2-1} \nabla w_k  \right)   + \int_{\partial\Omega} \dfrac{\partial \frac 1 2 |w|^2}{\partial \nu} |w|^{p-2} d\sigma\\
&=&\displaystyle   - \dfrac{4(p-1)}{p^2}\int_{\Omega} \left|  \nabla |w_k|^{\frac p 2}   \right|^2  + \frac 1 2\int_{\partial\Omega} \dfrac{\partial  |w|^2}{\partial \nu} |w|^{p-2} d\sigma.
\end{array}\end{equation*}
Hence, we obtain like in \cite{KS1}, proof of Theorem 1,
\begin{equation}\label{presk2}  -\varepsilon \sum_{k=1}^d  \int_\Omega |w_k|^{p-2}w_k\Delta w_k \geqslant \dfrac{4\varepsilon(p-1)}{p^2}\int_\Omega \left|  \nabla |\nabla u_\varepsilon |^{p/2}  \right|^2 - \dfrac{\varepsilon}{2}\int_{\partial\Omega} \dfrac{\partial |\nabla u_\varepsilon|^2}{\partial \nu}|\nabla u_\varepsilon|^{p-2} d\sigma,\end{equation}
and, combining \eqref{presk1} and \eqref{presk2},
\begin{multline}\label{Inter}
\dfrac{d \norme{\nabla u_\varepsilon}_p^p}{dt} + \dfrac{4\varepsilon(p-1)}{p} \int_\Omega \left|  \nabla |\nabla u_\varepsilon |^{p/2}  \right|^2 - \frac{ p\varepsilon} 2\int_{\partial\Omega} \dfrac{\partial |\nabla u_\varepsilon|^2}{\partial \nu}|\nabla u_\varepsilon|^{p-2} d\sigma \leqslant\\ C_6 \mathcal{F}(t) \Big(1 + \norme{u_\varepsilon}_{\infty}^3(1+\log^+(\norme{u_\varepsilon}_{W^{1,p}}))\Big)\norme{\nabla u_\varepsilon }_p^p
\end{multline}
Using Lemma \ref{Conv}, we have $\displaystyle - \frac{ p\varepsilon} 2 \int_{\partial\Omega} \dfrac{\partial |\nabla u_\varepsilon|^2}{\partial \nu}|\nabla u_\varepsilon|^{p-2} d\sigma \geqslant 0$, and thus, using the continuous embedding of $W^{1,p}(\Omega)$ into $L^\infty(\Omega)$, there exists a constant $C_8>0$ independent from $\varepsilon$ such that
\begin{equation}\label{Brique2}
\dfrac{d }{dt} \norme{ \nabla u_\varepsilon}_p^p \leqslant C_8\mathcal F(t) \left(1 + \norme{u_\varepsilon}_{W^{1,p}}^3(1+\log^+(\norme{ u_\varepsilon}_{W^{1,p}}))\right)\norme{\nabla u_\varepsilon}_p^p.
\end{equation}
We combine the bounds \eqref{Brique1} and \eqref{Brique2} to obtain, for some constant $C_9>0$ independent from $\varepsilon$,
\begin{equation}
\dfrac{d }{dt} \norme{u_\varepsilon}_{W^{1,p}}^p \leqslant C_9\mathcal F(t) \left(1 +  \norme{u_\varepsilon}_{W^{1,p}}^3(1+\log^+(\norme{u_\varepsilon}_{W^{1,p}}^p))\right)\norme{u_\varepsilon}_{W^{1,p}}^p.
\end{equation}
Hence, $\norme{u_\varepsilon}_{W^{1,p}}$ is a sub-solution to the ordinary differential equation
\begin{equation}\label{Some_ODE} w'(t) =  C_9\mathcal F(t) \Big(1 + w(t)^{\frac 3 p} \left(1 + \log^+(w(t)) \right) \Big) w(t).   \end{equation}
There is no term depending on $\varepsilon$ in the previous ODE. Hence, by Lemma \ref{upper_sol} applied with initial datum $\norme{u_0}_{W^{1,p}}^p$ there exists a time $T>0$ independent from $\varepsilon$ and there exists a constant $M_2(T)>0$ independent from $\varepsilon$ such that the following uniform estimate holds for all $t\in[0,T]$:
\[ \forall \varepsilon\in\R_+^*, \ \norme{u_\varepsilon(t)}_{W^{1,p}}\leqslant  M_2(T).\]
Using the continuous embedding of $W^{1,p}(\Omega)$ into $L^\infty(\Omega)$, we obtain the $L^\infty$ $\varepsilon-$uniform estimate: there exists $M_1(T)>0$ such that for all $t\in[0,T]$,
\[ \forall \varepsilon\in\R_+^*, \ \ \norme{u_\varepsilon(t)}_\infty < M_1(T).\]\qed
\begin{lemma}\label{lemme3}
	Assume $\Omega$ is bounded and convex, with a smooth boundary. Let $\varepsilon_0>0$. Let $u_0^\varepsilon,n_0^\varepsilon$ satisfy Hypothesis \ref{H0eps}. There exist two positive constants $M_3=M_3(T)$ and $M_4=M_4(T)$ depending on $T$ such that for all $\varepsilon\in(0,\varepsilon_0)$, \begin{equation} \label{L2B3} \sup_{t\in[0,T]}\norme{\dfrac{\partial u_\varepsilon}{\partial t}(t)}_{W^{-1,p}} < M_3(T)\qquad \mathrm{and} \qquad \sup_{t\in[0,T]}\norme{\dfrac{\partial n_\varepsilon}{\partial t}(t)}_{L^p} < M_4(T).\end{equation}
\end{lemma}

\noindent{\bf Proof.} First, the equation for $u$ reads
\[\dfrac{\partial u_{\varepsilon}}{\partial t} = \varepsilon\nabla \cdot \nabla u_{\varepsilon} - \nabla \cdot (u_{\varepsilon}\nabla  c_{\varepsilon} ) + g(u_{\varepsilon})n_{\varepsilon}u_{\varepsilon} - b(n_{\varepsilon})u_{\varepsilon}\]
Applying Lemma \ref{lemme2}, there exists $T>0$, $M_1(T)>0$ and $M_2(T)>0$ independent from $\varepsilon$ such that for all test function $\phi\in W^{1,q}(\Omega)$, with $\frac{1}{p}+\frac{1}{q}=1$, we have the following uniform bounds, independently of $\varepsilon$ and $t$,
\[|\pscal{g(u_{\varepsilon})n_{\varepsilon}u_{\varepsilon} - b(n_{\varepsilon})u_{\varepsilon}}{\phi}_{W^{-1,p},W^{1,q}}| \leqslant\displaystyle (G_0 \norme{n_0}_\infty + B_0) M_1(T) \norme{\phi}_{W^{1,q}} ;\]
\[
\begin{array}{rcl}
\displaystyle |\pscal{\nabla \cdot (u_{\varepsilon}\nabla c_{\varepsilon} )}{\phi}_{W^{-1,p},W^{1,q}}|&\leqslant&\displaystyle \int_{\Omega} u_{\varepsilon} | \nabla c_{\varepsilon} \cdot \nabla \phi |\\&\leqslant&\displaystyle M_1(T)   \norme{\nabla c_{\varepsilon}}_p \norme{\nabla \phi}_q \\&\leqslant&\displaystyle M_1(T)^2 K_\infty  \norme{\phi}_{W^{1,q}} ;
\end{array}
\]
\[
\begin{array}{rcl}
|\varepsilon\pscal{\nabla \cdot \nabla u_{\varepsilon}}{\phi}_{W^{-1,p},W^{1,q}}|&\leqslant&\displaystyle \varepsilon_0\int_{\Omega} | \nabla u_{\varepsilon} \cdot \nabla \phi | \\ &\leqslant&\displaystyle \varepsilon_0\norme{\nabla u_{\varepsilon}}_p \norme{\nabla \phi}_q \\ &\leqslant&\displaystyle \varepsilon_0M_2(T)\norme{\phi}_{W^{1,q}} ;\\
\end{array}
\]
where we used integration by parts, Hölder inequality, elliptic regularity (Lemma \ref{ElliptReg1}), boundary conditions and Lemma \ref{lemme2}. The result follows.\\
Second, the equation for $n$ reads 
\[  \dfrac{\partial n_{\varepsilon}}{\partial t} -\Delta n_\varepsilon =  - \gamma g(u_\varepsilon)n_\varepsilon u_\varepsilon.  \]
%From \eqref{NightWork} in the proof of Lemma \ref{lemme2}, we know that there exists a constant $C=C(\norme{n_0}_{W^{1,+\infty}})>0$ such that
%\[  \norme{\nabla u}_{W^{1,\infty}} \leqslant C (1+\norme{u}_\infty) (1+ t + 2\sqrt{t}) \leqslant C(1+M_1(T))(1+T+2\sqrt{T}). \]
We deduce from \cite{C2}, Section 5, equation (5.4) that, defining $Q_T=\Omega \times(0,T)$,
\[  \norme{n_\varepsilon}_{W^{2,1}_p}   \leqslant \lambda_p \left( \norme{u_0^\varepsilon}_{W^{2,p}} + \norme{\gamma g(u)nu}_{L^p(Q_T)}  \right)  ,  \]
for some constant $\lambda_p>0$. Using Hypothesis \eqref{H0eps}, we conclude that
\begin{equation}\label{Changfu}  \norme{n_\varepsilon}_{W^{2,1}_p}   \leqslant \lambda_p \left( \norme{u_0}_{W^{2,p}} + \gamma G_0\norme{n_0}_\infty M_1(T)  |\Omega|^\frac 1p T^\frac1p   \right)  ,  \end{equation}
Hence the result.
\qed

%
%Using theorem \ref{Conv}, we have $ - C_1\varepsilon \int_{\partial\Omega} \dfrac{\partial |\nabla u|^2}{\partial \mu}|\nabla u|^{p-2} d\sigma \geqslant 0$, and thus, using that $\norme{u}_\infty < M$, there exists a constant $C_3>0$ such that
%\begin{equation}
%\dfrac{d \norme{\nabla u}_p^p}{dt}  \leqslant C_3 \left(1 + log^+(\norme{\nabla u}^p_p)\right)\norme{\nabla u }_p^p
%\end{equation}
%Then as there is no term depending on $\varepsilon$ in the previous ODE, Gronwäll's lemma gives the uniform estimate
%\[\norme{\nabla u_\varepsilon(t)}_p\leqslant \e^{log(1+\norme{\nabla u_0}_p)\e^{Ct}},\]
%for some constant $C$ independent from $\varepsilon$.

\section{Solutions for the hyperbolic-elliptic-parabolic system}

In the following, we denote by $\mathscr{C}^\infty_0(\bar \Omega\times [0,T))$ the space of $\mathscr{C}^\infty$ functions on $\bar \Omega\times [0,T)$ whose support is a compact subset of $\bar \Omega\times [0,T)$.\\

We first define a notion of solution for problem \eqref{Intermediaire}:

\begin{definition}\label{TheDef}
	A strong $L^p$ solution of \eqref{Intermediaire} on $[0,T^*)$ is a triplet $(u,c,n)$ of functions such that there exists $p\in(1,+\infty)$, such that for all $T\in(0,T^*)$,
	\[ u\in  L^p(0,T;W^{1,p}(\Omega))\cap L^\infty(\Omega\times(0,T)), \quad c,n\in  L^p(0,T;W^{2,p}(\Omega))\cap L^\infty(\Omega\times(0,T)),  \]
	\[ \frac {\partial u}{\partial t},\ \frac {\partial n}{\partial t}\in  L^p(0,T;L^p(\Omega)), \qquad \nabla \cdot \Big( u\nabla  c \Big) \in  L^p(0,T;L^p(\Omega)),\]
	and $(u,c,n)$ solves \eqref{Intermediaire} in $L^p(0,T;L^p(\Omega))$ (and thus almost everywhere) with the homogeneous Neumann boundary conditions on $c$ and $n$ satisfied in the sense of $L^p(0,T;W^{1-\frac 1 p, p}(\partial\Omega))$.
\end{definition}

\begin{remark}\label{HandyRem}
	If for all $T\in(0,T^*)$, $(u,c,n)$ satisfies 
	\[ u\in  L^p(0,T;W^{1,p}(\Omega))\cap L^\infty(\Omega\times(0,T)), \quad c,n\in  L^p(0,T;W^{2,p}(\Omega))\cap L^\infty(\Omega\times(0,T)),  \]
	and if for all $\phi_1,\phi_2,\phi_3\in\mathscr{C}^\infty_0(\bar\Omega\times[0,T))$,
	\begin{equation}
	\int_{(0,T)\times\Omega} u \dfrac{\partial \phi_1}{\partial t} + \int_{(0,T)\times\Omega} u \nabla c \cdot \nabla \phi_1\\ + \int_{(0,T)\times\Omega} \left( g(u)n - b(n)\right)u\phi_1 = - \int_\Omega u_0\phi_1(\cdot,0),
	\end{equation}
	\begin{equation}
	\int_{(0,T)\times\Omega} \nabla  c  \cdot \nabla \phi_2 + \beta \int_{(0,T)\times\Omega} c\phi_2 = \alpha \int_{(0,T)\times\Omega} u \phi_2
	\end{equation}
	and
	\begin{equation}
	\int_{(0,T)\times\Omega} n \dfrac{\partial \phi_3}{\partial t} - \int_{(0,T)\times\Omega} \nabla n \cdot \nabla \phi_3  - \gamma \int_{(0,T)\times\Omega} g(u)nu\phi_3 = - \int_\Omega n_0\phi_3(\cdot,0).
	\end{equation}
	then by standard integration arguments, $(u,c,n)$ is a solution of \eqref{Intermediaire} on $[0,T^*)$ in the sense of definition \ref{TheDef}.
\end{remark}

\subsection{Existence of local-in-time solutions}

We are going to construct the solutions as the limit of a sequence of solutions of the approximate system \eqref{MimuraEps} via Aubin-Lions compactness theorem together with the uniform \textit{a priori} estimates in the lemmas \ref{lemme2} and \ref{lemme3}.

\begin{theorem}
	Assume $\Omega$ is bounded and convex, with a smooth boundary. Let $u_0,n_0$ satisfy Hypothesis \ref{H0}. There exists a unique non-negative solution $(u,c,n)$ of \eqref{Intermediaire} in the sense of Definition \ref{TheDef} on $[0,T]$.
\end{theorem}

\noindent{\bf Proof.}
\textit{Existence:\\}
Fix $p>\max(d,4)$. Let $(\varepsilon_k)_{k\in\N}$ be a decreasing sequence of positive numbers with limit 0. For clarity, we will also denote by $(\varepsilon_k)_{k\in\N}$ all its subsequences.\\
Fix $\eta\in(0,1)$. Let $(u_0^{\varepsilon_k})_{k\in\N}$ be a sequence of non-negative functions in $\mathscr{C}^{2,\eta}(\bar \Omega)$, satisfying homogeneous Neumann boundary conditions, which converges to $u_0$ in $W^{2,p}(\Omega)$. Let $(n_0^{\varepsilon_k})_{k\in\N}$ be a sequence of non-negative functions in $\mathscr{C}^{2,\eta}(\bar \Omega)$ such that $0\leqslant n_0^{\varepsilon_k}\leqslant \norme{n_0}_\infty$, satisfying homogeneous Neumann boundary conditions, which converges to $n_0$ in $W^{2,p}(\Omega)$.
We denote by $(u_{\varepsilon_k},c_{\varepsilon_k},n_{\varepsilon_k})$ the classical solution of system \eqref{MimuraEps} with diffusion parameter $\varepsilon = \varepsilon_k$ and initial condition $(u_0^{\varepsilon_k},n_0^{\varepsilon_k})$; this unique solution is given by Theorem \ref{thm:localE}.

By Lemma \ref{lemme2}, the sequence $(u_{\varepsilon_k})_{k\in\N}$ is uniformly bounded in $L^\infty(0,T;L^\infty(\Omega))$ by a constant $M_1(T)$ independent from $\varepsilon_k$ and thus bounded in $ L^p(0,T;L^p(\Omega))$. Also by Lemma \ref{lemme2}, the sequence $(\nabla u_{\varepsilon_k})_{k\in\N}$ is bounded in $L^\infty(0,T;L^p(\Omega))$  by a constant $M_2(T)$ independent from $\varepsilon_k$, and thus $(\nabla u_{\varepsilon_k})_{k\in\N}$ is bounded in $ L^p(0,T;L^p(\Omega))$.\\
The sequence $(u_{\varepsilon_k})_{k\in\N}$ is then bounded in $ L^p(0,T;W^{1,p}(\Omega))$, which is a reflexive space. Then, there exists $u\in  L^p(0,T;W^{1,p}(\Omega))$ such that, up to extraction of a subsequence
\begin{equation}\label{arg1}
u_{\varepsilon_k} \rightharpoonup u \qquad \mathrm{weakly \ in \ }  L^p(0,T;W^{1,p}(\Omega)).
\end{equation}
Moreover, by Lemma \ref{lemme3}, the sequence $(\dfrac{\partial u_{\varepsilon_k}}{\partial t})_{k\in\N}$ is bounded in $L^\infty(0,T;W^{-1,p}(\Omega))$ by a constant $M_3(T)$ independent from $\varepsilon_k$ and thus bounded uniformly in $L^2(0,T;W^{-1,p}(\Omega))$.
By Aubin-Lions compactness theorem (Lemma \ref{Aubin-Lions}) applied to $W^{1,p}(\Omega)$, $L^p(\Omega)$ and $W^{-1,p}(\Omega)$, the subspace
\[ \mathcal{Y} = \left\{ u\in  L^p(0,T;W^{1,p}(\Omega)) \ | \  	\dfrac{\partial u}{\partial t}\in L^2(0,T;W^{-1,p}(\Omega)) \right\} \]
is compactly embedded in $ L^p(0,T;L^p(\Omega))$. Therefore, up to extraction of a subsequence, $(u_{\varepsilon_k})_{k\in\N}$ converges strongly in $ L^p(0,T;L^p(\Omega))$ towards $u$.

Since $p>2$, the sequence $(u_{\varepsilon_k})_{k\in\N}$ converges towards $u$ in $L^2(Q_T)$, where $Q_T=\Omega\times (0,T)$. Let $c$ be the classical solution
\[  -\Delta c + \beta c = \alpha u  \quad \mathrm{in}\ \Omega, \qquad \dfrac{\partial c}{\partial \nu}=0 \quad \mathrm{on}\ \partial\Omega .\]
Denote $U_k=u-u_{\varepsilon_k}$ and $C_k= c-c_{\varepsilon_k}$. We have
\[  -\Delta c_{\varepsilon_k}   + \beta c_{\varepsilon_k}  = \alpha u_{\varepsilon_k},    \]
so that,
\[    -\Delta C_k  + \beta C_k = \alpha U_k  \]
We multiply this equation by $C_k$ and we integrate in space and time. We obtain
\[  \int_0^T \int_{\Omega}  \left| \nabla C_k \right|^2  + \beta\int_0^T \int_{\Omega} C_k^2  = \alpha \int_0^T \int_{\Omega} U_kC_k  \]
and then
\[  \int_0^T \int_{\Omega}  \left| \nabla C_k \right|^2  + \frac\beta2\int_0^T \int_{\Omega} C_k^2  \leqslant \frac{\alpha^2}{2\beta} \int_0^T \int_{\Omega} U_k^2.  \]
Since the right-hand side converges to 0 when $k$ tends to $+\infty$, the sequence $(c_{\varepsilon_k})_{k\in\N}$ converges towards $c$ in $L^2(0,T;H^1(\Omega))$.

Using \eqref{Changfu}, we deduce that $n$ is bounded uniformly in $\varepsilon$ in $L^\infty(0,T;W^{2,p}(\Omega))$ and thus in $L^p(0,T;W^{2,p}(\Omega))$. Then, there exists $n\in  L^p(0,T;W^{2,p}(\Omega))$ such that, up to extraction of a subsequence
\begin{equation}\label{arg2}
	n_{\varepsilon_k} \rightharpoonup n \qquad \mathrm{weakly \ in \ }  L^p(0,T;W^{2,p}(\Omega)).
\end{equation} 
By Lemma \ref{lemme3}, the sequence $(\dfrac{\partial n_{\varepsilon_k}}{\partial t})_{k\in\N}$ is bounded in $L^\infty(0,T;L^p(\Omega))$ by a constant $M_4(T)$ independent from $\varepsilon_k$ and thus bounded uniformly in $L^2(0,T;L^p(\Omega))$. Applying the Aubin-Lions compactness theorem (Lemma \ref{Aubin-Lions}), we deduce that $n_{\varepsilon_k}$ converges strongly towards $n$ in $L^p(0,T;W^{1,p}(\Omega))$

Since classical solutions are also weak solutions, for all $\phi_1,\phi_2,\phi_3\in\mathscr{C}^\infty_0(\bar\Omega\times[0,T))$, for all $k\in\N$ (after all required extractions),
\begin{multline}\label{lim1}
\int_{(0,T)\times\Omega} u_{\varepsilon_k} \dfrac{\partial \phi_1}{\partial t} - \varepsilon_k \int_{(0,T)\times\Omega} \nabla u_{\varepsilon_k} \cdot \nabla \phi_1  + \int_{(0,T)\times\Omega} u_{\varepsilon_k} \nabla c_{\varepsilon_k} \cdot \nabla \phi_1\\ + \int_{(0,T)\times\Omega} \left( g(u_{\varepsilon_k})n_{\varepsilon_k} - b(n_{\varepsilon_k})\right)u_{\varepsilon_k}\phi_1 = - \int_\Omega u_0^{\varepsilon_k}\phi_1(\cdot,0),
\end{multline}
\begin{equation}\label{lim_2}
\int_{(0,T)\times\Omega} \nabla c_{\varepsilon_k} \cdot \nabla \phi_2 + \beta\int_{(0,T)\times\Omega} c_{\varepsilon_k} \phi = \alpha \int_{(0,T)\times\Omega} u_{\varepsilon_k} \phi_2
\end{equation}
and
\begin{equation}\label{lim3}
\int_{(0,T)\times\Omega} n_{\varepsilon_k} \dfrac{\partial \phi_3}{\partial t} - \int_{(0,T)\times\Omega} \nabla n_{\varepsilon_k} \cdot \nabla \phi_3  - \gamma \int_{(0,T)\times\Omega} g(u_{\varepsilon_k})n_{\varepsilon_k}u_{\varepsilon_k}\phi_3 = - \int_\Omega n_0^{\varepsilon_k}\phi_3(\cdot,0).
\end{equation}
Since $g$ and $b$ are continuously differentiable, we have by the strong convergence of $u_{\varepsilon_k}$ and $n_{\varepsilon_k}$ in $L^p(0,T;L^p(\Omega))$ and by the Lebesgue dominated convergence theorem that $g(u_{\varepsilon_k})n_{\varepsilon_k}u_{\varepsilon_k}$ converges towards $g(u)nu$ and $b(n_{\varepsilon_k})u_{\varepsilon_k}$ converges towards $b(n)u$ almost everywhere in $\Omega\times [0,T]$. Using also the previously obtained weak convergences, we pass to the limit in \eqref{lim1}, \eqref{lim_2} and \eqref{lim3} and we obtain
\begin{equation}
\int_{(0,T)\times\Omega} u \dfrac{\partial \phi_1}{\partial t} + \int_{(0,T)\times\Omega} u \nabla c \cdot \nabla \phi_1\\ + \int_{(0,T)\times\Omega} \left( g(u)n - b(n)\right)u\phi_1 = - \int_\Omega u_0\phi_1(\cdot,0),
\end{equation}
\begin{equation}
\int_{(0,T)\times\Omega} \nabla  c  \cdot \nabla \phi_2 + \beta\int_{(0,T)\times\Omega} c\phi_2 = \alpha\int_{(0,T)\times\Omega} u \phi_2
\end{equation}
and
\begin{equation}
\int_{(0,T)\times\Omega} n \dfrac{\partial \phi_3}{\partial t} - \int_{(0,T)\times\Omega} \nabla n \cdot \nabla \phi_3  - \gamma \int_{(0,T)\times\Omega} g(u)nu\phi_3 = - \int_\Omega n_0\phi_3(\cdot,0).
\end{equation}
The functions $(u,c,n)$ are weak or strong limits of sequences of non-negative functions, they are thus non-negative as well. Moreover, since $n$ is the strong limit in $ L^p(0,T;L^p(\Omega))$ of a sequence of functions satisfying $0<  n_{\varepsilon_k} \leqslant \norme{n_0}_\infty$, 
then passing to the limit $k\to+\infty$, we obtain
\[ 0 \leqslant n \leqslant \norme{n_0}_\infty.  \]
Using the regularity properties of $u$, $c$ and $n$ and in view of Remark \ref{HandyRem}, we deduce that  $(u,c,n)$ is a solution of \eqref{Intermediaire} in the sense of Definition \ref{TheDef}.\vspace{0.3cm}\\ 
\noindent\textit{Uniqueness:\\}
Assume there exist two solutions $(u_1,c_1,n_1)$ and $(u_2,c_2,n_2)$ on $\Omega\times[0,T]$. Denote $\bar u = u_1 - u_2$, $\bar c = c_1 - c_2$ and $\bar n = n_1 - n_2$. 
By definition \ref{TheDef}, we have
\[  \dfrac{\partial \bar u}{\partial t},\dfrac{\partial \bar n}{\partial t}\in  L^p(0,T;L^p(\Omega)).   \]
Thus, the equations for $(u_1,c_1,n_1)$ and $(u_2,c_2,n_2)$ are satisfied in $L^p(0,T;L^p(\Omega))$ and hence satisfied almost everywhere in $\Omega \times [0,T]$. Henceforth, we implicitly consider these equations on the subset of $[0,T]$ whereon they are defined. The triplet $(\bar u, \bar c, \bar n)$ is a solution of
\begin{equation}\label{eq1f} \dfrac{\partial \bar u}{\partial t} = - \nabla \cdot (\bar u \nabla c_1) - \nabla \cdot \Big( u_2 \nabla  \bar c \Big) + g(u_1)n_1u_1 - g(u_2)n_2u_2 - b(n_1)u_1 + b(n_2)u_2, \end{equation}
\begin{equation}\label{eq2f} 0 = \Delta \bar c - \beta \bar c + \alpha \bar u, \end{equation}
\begin{equation}\label{eq3f} \dfrac{\partial \bar n}{\partial t} = \Delta \bar n - \gamma \Big(g(u_1)n_1u_1 - g(u_2)n_2u_2\Big).  \end{equation}
%Multiplying \eqref{eq2f} by $\bar c$ and integrating over $\Omega$, we obtain
%\[ \norme{\nabla \bar c}^2_2 + \norme{ \bar c}^2_2 \leqslant \norme{ \bar u}_2\norme{ \bar c}_2 \leqslant \dfrac{1}{2}\norme{ \bar u}^2_2 + \dfrac{1}{2}\norme{ \bar c}^2_2, \]
%and that implies
%\begin{equation}\label{Help} \norme{\nabla \bar c}^2_2+\dfrac{1}{2}\norme{ \bar c}^2_2 \leqslant \dfrac{1}{2}\norme{ \bar u}^2_2. \end{equation} 
First, we have
	\begin{equation}\label{gbound}
	|g(u_1)u_1n_1-g(u_2)u_2n_2| \leqslant 
	G_0 \norme{u_1}_\infty|\bar n |  +
	\sup_{s\in[0,T]}\norme{n_2(s)}_\infty 
	 (\norme{g'}_\infty\norme{u_2}_\infty+G_0)|\bar u|,
	\end{equation}
	and
	\begin{equation}\label{bbound}
	|b(n_1)u_1-b(n_2)u_2|\leqslant   
	B_0|\bar u|
	+\norme{u_2}_\infty \norme{b'}_\infty |\bar n|.
	\end{equation}
	We multiply \eqref{eq3f} by $\bar n$ and we integrate. There exists $\kappa_3>0$ such that,
	\begin{equation}\label{IneqN2}
	\dfrac{d}{dt}\int_\Omega \bar n^2 + 2\int_{\Omega} |\nabla \bar n|^2 \leqslant \kappa_3 \left(  \int_\Omega \bar u^2  + \int_\Omega \bar n^2 \right).
	\end{equation}
We are going to multiply \eqref{eq1f} by $\bar u$ and to integrate. Hence, we make the following preliminary computations.
Using Lemma \ref{ElliptReg1}, we have
\begin{equation}\label{FirstBound}
 - \int_\Omega \nabla \cdot(\bar u \nabla c_1)\bar u =  -\dfrac{1}{2} \int_\Omega |\bar u|^2 \Delta c_1 
  = -\dfrac{1}{2} \int_\Omega |\bar u|^2 ( \beta c_1 - \alpha u_1 )  \leqslant \dfrac{1}{2}  (\alpha+K_\infty\beta)\norme{u_1}_\infty \norme{\bar u}_2^2   .
\end{equation}
We also compute
\begin{equation}\label{FirstAndHalfBound}
-\int_\Omega \nabla \cdot \Big( u_2 \nabla  \bar c \Big)\bar u = -\int_\Omega \bar u \nabla u_2 \cdot \nabla  \bar c   - \int_\Omega  u_2 \bar u \Delta \bar c  
\leqslant  \int_\Omega |\nabla u_2| |\nabla \bar c||\bar u| + \norme{u_2}_\infty \int_\Omega |\Delta \bar c||\bar u| .
\end{equation}
We need a bound for each term of the last sum. We deduce from Lemma \ref{ElliptReg1},
\begin{multline}\label{SecondBound}
\int_\Omega |\Delta \bar c ||\bar u| = \int_\Omega | \alpha\bar u - \beta\bar c ||\bar u| \leqslant \alpha\norme{\bar u}_2^2 + \frac{\beta^2}2 \norme{\bar c}_2^2 + \frac12\norme{\bar u}_2^2 \leqslant (\alpha+\frac 12 + \frac{\beta^2K_2}2 )\norme{u}_2^2,
\end{multline}
For the first term in \eqref{FirstAndHalfBound}, we are going to use different bounds depending on the space dimension $d$.

\noindent\textit{Case $d\geqslant 3$:}
Using the extended Hölder inequality (cf. \cite{B2}, Remark 2, p. 93), the Sobolev embedding of $H^1(\Omega)$ into $L^\frac{2d}{d-2}(\Omega)$ (cf. \cite{B2}, Corollary 9.14) and Lemma \ref{ElliptReg1}, we have for some large enough constants $C_1,C_2>0$,
\begin{equation}\label{ThirdBound}
\int_\Omega |\nabla u_2| |\nabla \bar c||\bar u|  \leqslant \norme{\nabla u_2}_d  \norme{\nabla \bar c}_{\frac{2d}{d-2}} \norme{\bar u}_2 \leqslant   C_1 \norme{\nabla u_2}_d \norme{\nabla \bar c}_{H^1}  \norme{\bar u}_2 \leqslant  C_2 \norme{\nabla u_2}_d   \norme{\bar u}_2^2.
\end{equation}

The quantities $\norme{u_1}_\infty, \norme{u_2}_\infty, \norme{n_2}_\infty, \norme{n_2}_\infty, \norme{\nabla u_2}_d $ are uniformly bounded on $[0,T]$. We multiply equation \eqref{eq1f} by $\bar u$ and we integrate over $\Omega$. Then, using the inequalities \eqref{gbound}, \eqref{bbound}, \eqref{FirstBound}, \eqref{FirstAndHalfBound}, \eqref{SecondBound} and \eqref{ThirdBound}, we deduce that there exists a constant $\kappa_4>0$ large enough such that 
\[  \dfrac{d }{d t}\norme{\bar u}_2^2 \leqslant \kappa_4\left( \norme{\bar u}^2_2 + \norme{\bar n}^2_2\right). \]
Adding inequality \eqref{IneqN2}, we deduce that there exists a constant $\kappa_5>0$ such that
\[  \dfrac{d }{d t}\left(\norme{\bar u}_2^2 + \norme{\bar n}^2_2 \right) \leqslant \kappa_5\left( \norme{\bar u}^2_2 + \norme{\bar n}^2_2\right). \]

Since $\norme{\bar u(0)}_2^2  + \norme{\bar n(0)}^2_2=0$, Gronwall's lemma implies that for almost every time $t$ in $[0,T]$, $\bar u(t) = 0$ and $\bar n(t) =0$ in $L^p(\Omega)$. Thus, we have $\bar u = \bar n = 0$ in $ L^p(0,T;L^p(\Omega))$.

\noindent\textit{Case $d = 1$ or $d=2$:}
Using the Sobolev embedding of $H^1(\Omega)$ into $L^4(\Omega)$ (cf. \cite{B2}, Corollary 9.14) and Lemma \ref{ElliptReg1}, we have for some constants $C_4,C_5>0$
\begin{equation}\label{ThirdBoundBis}
\displaystyle \int_\Omega |\nabla u_2|  |\nabla \bar c| |\bar u| \leqslant\displaystyle   \norme{\nabla u_2}_4  \norme{\nabla \bar c}_4 \norme{\bar u}_2\\
\leqslant C_4 \norme{\nabla u_2}_4 \norme{\nabla \bar c}_{H^1}  \norme{\bar u}_2\\
\leqslant C_5 \norme{\nabla u_2}_4   \norme{\bar u}_2^2.
 \end{equation}
Using the same reasoning as above but replacing \eqref{ThirdBound} by \eqref{ThirdBoundBis}, we prove $\bar u = \bar n= 0$ in the space $ L^p(0,T;L^p(\Omega))$.

Last, in all dimensions, $\bar u = 0$. The function $\bar c$ is then solution of
\[  -\Delta \bar c + \beta \bar c = 0,\]
with homogeneous Neumann boundary condition. Therefore, $ \bar c = 0$ in the space $ L^p(0,T;L^p(\Omega))$, which completes the proof of uniqueness.
\qed

%%%%%%%%%%%%%%%%%%%%%%%%%%%%%%%%%%%%%%%%%%%%%%%%%%%%%%%%%%%%%%%

\section{Blow-up in finite time in the hyperbolic case}
\label{sec: blowup}

We focus here on finite time blow-up of the solutions of \eqref{Intermediaire}. We build on the methods of Winkler \cite{W1}, Lankeit \cite{L3} and Kang and Stevens \cite{KS1}.

\begin{theorem}
	Assume $\Omega$ is bounded with a smooth boundary. Let $p\in(1,+\infty)$. Let $(u,c,n)$ be a strong $L^p$ solution of \eqref{Intermediaire} in the sense of Definition \ref{TheDef}. Then there exists a constant $C_p=C_p(\norme{u_0}_1,\norme{n_0}_1)>0$ depending on $\norme{u_0}_1$, $\norme{n_0}_1$, $p$ and $\Omega$ such that if 
	\[ \norme{u_0}_p > C_p(\norme{u_0}_1,\norme{n_0}_1), \] 
	then the solution $u$ blows-up in finite time in $L^p$ and $L^\infty$:
	\[ \lim_{t\to T_{\mathrm{max}}}\norme{u(t)}_p =  \lim_{t\to T_{\mathrm{max}}}\norme{u(t)}_\infty = +\infty. \]
\end{theorem}

{\bf Proof.}
By Definition \ref{TheDef}, $\dfrac{\partial u}{\partial t}\in L^p(0,T;L^p{\Omega})$. Thus, we can multiply the first equation of \eqref{Intermediaire} with the $L^{\frac{p}{p-1}}$ function $u^{p-1}$ and compute
%\[    \pscal{ \dfrac{\partial u}{\partial t}}{u^{p-1}}_{ W^{-1,p},W^{1,\frac{p}{p-1}} }  = - \int_\Omega\nabla \cdot (u\nabla c)u^{p-1}  +\int_\Omega (g(u)n-b(n))u^p .\]
%Then, after integration by parts and simplifications, we obtain
\[ \dfrac{1}{p}\dfrac{d }{d t}\int_\Omega u^p = -  \dfrac{p-1}{p}\int_\Omega u^p\Delta c + \int_\Omega (g(u)n-b(n))u^p.\]
Using the second equation of \eqref{Intermediaire}, we have
\begin{equation} \dfrac{1}{p}\dfrac{d }{d t}\int_\Omega u^p = \alpha \dfrac{p-1}{p}\int_\Omega u^{p+1} - \beta \dfrac{p-1}{p}\int_\Omega cu^p + \int_\Omega (g(u)n-b(n))u^p.\label{zzz}\end{equation}

\noindent\textit{Case $d\geqslant 2$:\\}
We deduce from a Sobolev injection (\cite{B2}, p. 284, Corollary 9.14) and elliptic regularity (Lemma \ref{ElliptReg1}) that there exist constants $\tilde C_1,C_1>0$ such that
\begin{equation}\label{j1}   \norme{c}_{p+1}\leqslant \tilde C_1\norme{c}_{W^{1,q}} \leqslant C_1 \norme{u}_q,\end{equation}
where $\dfrac{1}{q} = \dfrac{1}{p+1}+\dfrac{1}{d}$, and by interpolation (\cite{B2}, p. 93, Remark 2),
\begin{equation}\label{j2} \norme{u}_q\leqslant \norme{u}_1^{1-\theta} \norme{u}^\theta_{p+1},\end{equation}
where $\dfrac{1}{q} = \dfrac{1-\theta}{1} + \dfrac{\theta}{p+1} $, that is to say $ \theta = \dfrac{p(d-1)-1}{pd} \in(0,1)$.\\
Moreover, with Hölder inequality applied with the exponents $p+1$ and $\dfrac{p+1}{p}$, we have, in view of \eqref{j1} and \eqref{j2},
\[ \int_{\Omega} cu^p \leqslant \left( \int_{\Omega} c^{p+1} \right)^\frac{1}{p+1} \left( \int_{\Omega} (u^p)^{\frac{p+1}{p}} \right)^{\frac{p}{p+1}}  = \norme{c}_{p+1}\norme{u}_{p+1}^p. \]
Applying these inequalities to \eqref{zzz}, we have
\[ \dfrac{1}{p}\dfrac{d }{d t} \norme{u}_p^p \geqslant \alpha \dfrac{p-1}{p} \norme{u}_{p+1}^{p+1} - \beta  C_1\dfrac{p-1}{p}\norme{u}_1^{1-\theta}\norme{u}_{p+1}^{p+\theta} - B_0 \norme{u}_p^p.  \]
Using the inequality $xy\leqslant \dfrac{x^a}{\varepsilon^\frac{a}{b} a}+\dfrac{\varepsilon y^b}{b}$ that holds for every $x,y,\varepsilon,a,b>0$ such that $\dfrac{1}{a}+\dfrac{1}{b}=1$, with the choice $x=\norme{u}_1^{1-\theta}$, $y=\norme{u}_{p+1}^{p+\theta}$, $a=\dfrac{p+1}{1-\theta}$ and $b=\dfrac{p+1}{p+\theta}$, we have
\begin{equation}
\dfrac{1}{p}\dfrac{d }{d t} \norme{u}_p^p \geqslant \alpha \dfrac{p-1}{2p} \norme{u}_{p+1}^{p+1} - \beta  C_{2}\norme{u}_1^{p+1} - B_0 \norme{u}_p^p, \label{yyy}
\end{equation}
where whe chose $\varepsilon$ such that $ \alpha\dfrac{p-1}{p} - \varepsilon\dfrac{p+\theta}{p+1}C_1 \beta \geqslant \alpha\dfrac{p-1}{2p} $ and defined $C_2 =\dfrac{p-1}{p} C_1\dfrac{1-\theta}{\varepsilon^\frac{p+\theta}{1-\theta}(p+1)}$.\\

\noindent\textit{Case $d=1$:\\}
We have again
\[ \int_{\Omega} cu^p \leqslant \left( \int_{\Omega} c^{p+1} \right)^\frac{1}{p+1} \left( \int_{\Omega} (u^p)^{\frac{p+1}{p}} \right)^{\frac{p}{p+1}}  = \norme{c}_{p+1}\norme{u}_{p+1}^p. \]
Using the inequality $xy\leqslant \dfrac{x^a}{\varepsilon^\frac{a}{b} a}+\dfrac{\varepsilon y^b}{b}$ with the choice $x=\norme{c}_{p+1}$, $y=\norme{u}_{p+1}^p$, $a=p+1$ and $b=\dfrac{p+1}{p}$, we have
\[ \norme{c}_{p+1}\norme{u}_{p+1}^p \leqslant  \dfrac{1}{\varepsilon^p(p+1)}\norme{c}_{p+1}^{p+1}  + \varepsilon\dfrac{p}{p+1}\norme{u}_{p+1}^{p+1}. \]
Let $\eta>0$. By Lemma \ref{ElliptReg2} (the hypothesis $u\in C(\bar \Omega)$ is satisfied thanks to the Sobolev injection of $W^{1,p}(\Omega)$ into $C(\bar \Omega)$ in the one dimensional case), there exists $C(\eta,p)$ such that
\[  \norme{c}_{p+1}^{p+1} \leqslant \eta \norme{u}_{p+1}^{p+1} + C(\eta,p)\norme{u}_1^{p+1}. \]
If we choose $\varepsilon$ and $\eta$ small enough, there exists a constant $\tilde C_2>0$ such that we deduce from \eqref{zzz} the estimate
\begin{equation}
\dfrac{1}{p}\dfrac{d }{d t} \norme{u}_p^p \geqslant \alpha \dfrac{p-1}{2p} \norme{u}_{p+1}^{p+1} - \beta  \tilde C_{2}\norme{u}_1^{p+1} - B_0 \norme{u}_p^p, \label{yyyBis}
\end{equation}
in the one dimensional case.\\

\noindent\textit{End of the proof for both $d=1$ and $d\geqslant 2$:\\}
Consider either \eqref{yyyBis} in the case $d=1$ or \eqref{yyy} in the case $d\geqslant 2$.\\
Since
\[ \dfrac{d }{d t}\int_\Omega \Big( u + \dfrac{1}{\gamma} n  \Big) = -\int_\Omega b(n)u \leqslant 0,\]
there exists a constant $C_{3}>0$ such that for all $t>0$
\[ \norme{u(t,\cdot)}_1 \leqslant C_{3}= C_3(\norme{u_0}_1,\norme{n_0}_1).\]
By continuous injection of $L^{p+1}$ into $L^p$, there exists a constant $C_4>0$ depending on $\Omega$ and $p$ such that
\[ \int_\Omega u^{p+1} \geqslant C_{4} \bigg( \int_\Omega u^p \bigg)^{\frac{p+1}{p}} .\]
The Inequalities \eqref{yyy} and \eqref{yyyBis} then become
\[ \dfrac{1}{p}\dfrac{d }{d t} \norme{u}_p^p \geqslant C_{5} \norme{u}_p^{p+1}-B_0\norme{u}_p^p-C_6,\]
where $C_5,C_6>0$ depend on $p$, $\alpha$, $\beta$, $\Omega$, $u_0$, $n_0$ and all previous constants, and this gives
\[ \dfrac{1}{p}\dfrac{d }{d t} \norme{u}_p^p \geqslant  \norme{u}_p^{p}\Big( C_{5} \norme{u}_p - B_0 \Big) - C_6.\]
Thus, if $\norme{u_0}_p$ is large enough, $\norme{u(t)}_p>\norme{u_0}_p$ for all $t>0$ and then $\norme{u}_p^p$ is a super-solution to an initial value problem of the form
\[  y'(t)= \alpha_1 y(t)^{1+\frac{1}{p}} - \alpha_2y(t) - \alpha_3, \]
where $\alpha_1,\alpha_2, \alpha_3>0$. Using Lemma \ref{BlowupODE}, we deduce that if $\norme{u_0}_p$ is large enough, there exists $T_{\mathrm{max}}>0$ such that
\[\lim_{t\to T_{\mathrm{max}}}\norme{u(\cdot,t)}_p=+\infty,\]
and, $\Omega$ being bounded,
\[\lim_{t\to T_{\mathrm{max}}}\norme{u(\cdot,t)}_\infty=+\infty.\]
%The blow-up rate is given by Lemma \ref{BlowupODE} and the fact that super-solutions to ODEs obviously have a higher blow-up rate than their sub-solution because their derivative is higher.
\qed

%%%%%%%%%%%%%%%%%%%%%%%%%%%%%%%%%%%%%%%%%%%%%%%%%%%%%%%%%%%%%%%%%

%
%\section{Global existence of solutions for the delayed model}
%
%
%\subsection{A criterion for the maximal time of existence}
%
%
%
%\subsection{Super-solutions and control over the firing rate}
%

\section{Conclusion and open problems}

This article focuses on a hyperbolic-elliptic-parabolic system of partial differential equations that represents pure chemotactic movement of self-organising E. Colies colonies. The original model of three parabolic equations \cite{AMM} proved to be adequate at describing the process by which stable geometrical colonial patterns form. It was not obvious that a hyperbolic version would be tractable. Our work is a first step towards the understanding of the pure chemotaxis regime. In order to carry out our analysis, we built on recent works on hyperbolic-elliptic chemotaxis systems \cite{W1}, \cite{L3}, \cite{KS1}, \cite{FGM}.

First, we constructed local-in-time solutions for a parabolic-elliptic-parabolic approximate system and we provided uniform \textit{a priori} estimates on these solutions. These estimates are uniform both in term of the diffusion coefficient of the first parabolic equation and in term of the initial conditions.

Then, we used these estimates and a standard compactness method involving Aubin-Lions theorem and we constructed solutions to the hyperbolic-elliptic-parabolic system. These solutions require mild assumptions on the parameters and functions of the system which correspond with the assumptions which were made in \cite{CHKMR}. 

Last, we studied the case of concentrated enough initial data and we proved that the corresponding solutions blow-up in finite time in $L^p$ and $L^\infty$ norms in all dimensions, including dimension 1.\\

Several questions regarding the solutions of this hyperbolic-elliptic-parabolic system remain open. We don't know yet if a non-linear chemotactic sensitivity function $\chi$ could prevent blow-up and we don't have results on the long time behaviour of global-in-time solutions. The asymptotic pattern could be easier to study than in the full parabolic case. Numerical simulations could also provide insights regarding the final pattern of global-in-time solutions.

Another interesting question is whether we can derive results on the parabolic-elliptic-parabolic system using the solutions of the hyperbolic-elliptic-parabolic system in the way of \cite{W1} and \cite{L3}. There are no theoretical obstacle as the solutions of the third equation are uniformly bounded and this could help to study the final pattern in the parabolic-elliptic-parabolic case and finally the full parabolic one, which is the initial goal of the model \cite{AMM}.

\appendix

\section{Classical results}

For the sake of clarity, we state here classical results that we use in our proofs.

\subsection{Elliptic estimates}

\begin{lemma}\label{ElliptReg1}
	Let $\Omega$ be a bounded open subset of $\R^d$ with a smooth boundary.
	\begin{itemize}
		\item Let $u\in L^{\infty}(\Omega)$. For all $p > 1$ and $\eta \in (0,1)$, there exists a unique solution $c \in W^{2,p}(\Omega)  \cap C^{1 + \eta}(\bar\Omega)$ of the boundary value problem
		\begin{equation}\label{A1}
		\left\{\begin{array}{rcll}
		-\Delta c +\beta c &=& \alpha u  &\quad \text{for}\quad x\in\Omega\\
		\dfrac{\partial c}{\partial \nu}&=&0 &\quad  \text{for}\quad x\in\partial\Omega,
		\end{array}\right.
		\end{equation}
		which is such that for all $p\in(1,+\infty)$,
		\begin{equation}\label{est2}
		\norme{c}_{W^{2,p}} \leqslant K_p\norme{u}_{L^{p}} \qquad \text{and} \qquad \norme{c}_ {C^{1 + \eta}} \leqslant K_\infty \norme{u}_\infty, 
		\end{equation}
		for some positive constants $K_p$, and $K_\infty$ depending on $\Omega$, $\alpha$, $\beta$, $d$ and $p$. If $u$ is nonnegative, then $c$ is nonnegative.
		\item If we assume that  $u\in W^{1,p}(\Omega)$, the solution $c$ of problem \eqref{A1} is such that
		\begin{equation}\label{est3}
		\norme{c}_{W^{3,p}} \leqslant K_h \norme{u}_{W^{1,p}}, 
		\end{equation}
		for some positive constant $K_h$ depending on $\Omega$, $\alpha$, $\beta$, $d$ and $p$.
		\item Let $\eta \in (0,1)$. There exists a constant $C>0$ depending on $\eta$, $\Omega$, $\alpha$, $d$ and $\beta$ such that for all $u\in C^{0,\eta}(\bar \Omega)$, the solution $c$ of problem \eqref{A1} is such that
		\begin{equation}\label{A4} \norme{c}_{C^{2,\eta}} \leqslant C \norme{u}_{C^{0,\eta}}.  \end{equation}
	\end{itemize}
\end{lemma}

{\bf Proof.} The existence and uniqueness of a solution $c \in {H^1}(\Omega)$ follows from applying Lax-Milgram theorem. The fact that $c \in W^{2,p}(\Omega)$ for all $p > 1$ can be deduced from  \cite{ADN}, whereas 
the last estimate in (\ref{est2}) follows from the embedding $W^{2,p}(\Omega)  \subset  C^{1 + \eta}(\bar\Omega)$ with $\eta = 1 - d/p$ if $p > d$.  The estimate  (\ref{est3}) also follows from \cite{ADN}.

The estimate \eqref{A4} follows from \cite{G}, Theorem 6.30. \qed

\bigskip

Next we define ${log^+} k = log~k$ for $k > 1$ and ${log^+} k = 0$ otherwise. One can find the following estimate in \cite{KS1}, \cite{OT}, Lemma 2.8, and \cite{Private}.

\begin{lemma}\label{BMO1}
	Let $\Omega$ be a bounded open set of $\R^d$ with a smooth boundary. Let $1<p<+\infty$ and $s>d/p$. There exists a constant $C$ depending on s, d and p such that for all $f\in W^{s,p}(\Omega)$,
	\[  \norme{f}_\infty \leqslant C\big(1+\norme{f}_{BMO}(1+\log^+(\norme{f}_{W^{s,p}}))\big).  \]
\end{lemma}   

\medskip

The following estimate is given by \cite{KS1} in the proof of Lemma 6. According to \cite{KS1} and \cite{Private2}, it can also be derived from \cite{S}, p 178 by choosing for the singular integral operator $T:L^\infty\to BMO$ the operator which maps $u$ to the Hessian of $c$.

\begin{lemma}[Estimates of singular integrals]\label{BMO2} Let $\Omega$ be a bounded open subset of $\R^d$ with a smooth boundary. Let $p\in(1,+\infty)$. Let $u\in L^\infty(\Omega)\cap W^{1,p}(\Omega)$. If the function $c$ is solution of
	\begin{equation}
	\left\{\begin{array}{rcll}
	-\Delta c +\beta c &=& \alpha u  &\quad \text{for}\quad x\in\Omega\\
	\dfrac{\partial c}{\partial \nu}&=&0 &\quad  \text{for}\quad x\in\partial\Omega
	\end{array}\right.
	\end{equation}
	with Dirichlet or Neumann boundary condition, then
	\[ \norme{\nabla^2 c}_{BMO} \leqslant C\norme{u}_\infty.\]
%	 \qquad \mathrm{and} \qquad  \norme{\nabla^2 c}_{W^{1,p}} \leqslant C\norme{\nabla u}_p
\end{lemma}

In the one dimensional case, we will use a more specific lemma. Its proof can be find in \cite{W1}, Lemma 2.2.

\begin{lemma}\label{ElliptReg2}
	Let $p>0$. Assume $d=1$ and $\Omega =\, (a,b)\,$. For all $\eta>0$ there exists $C(\eta,p)>0$ such that whenever $u\in C(\bar \Omega)$ is non-negative, the solution $c$ of 
	\begin{equation}
	\left\{\begin{array}{rcll}
	- c'' +\beta c &=& \alpha u  &\quad \text{for}\quad x\in (a,b)\\
	c'(a)=c'(b)&=&0 &
	\end{array}\right.
	\end{equation}
	satisfies
	\[  \int_\Omega c^{p+1} \leqslant \eta \int_\Omega u^{p+1} + C(\eta,p) \left(\int_\Omega u\right)^{p+1}.   \]
\end{lemma}

%
%
%\begin{theorem}[Agmon-Douglis-Nirenberg] \label{ElliptReg1}
%	Let $p\in(1,+\infty)$, $\alpha,\beta\in\R_+^*$ and $m\in\N$. If $\Omega$ is bounded and has a $\mathscr{C}^{m+2}$ boundary, if $u\in W^{m,p}(\Omega)$ there exists a unique classical solution $c\in W^{2,p}(\Omega)$ to
%	\[ -\Delta c + \beta c = \alpha u \]
%	with Dirichlet or Neumann boundary condition and there exists a constant $C\in\R_+^*$ independent from $u$ such that
%	\[ \norme{c}_{W^{m+2,p}} \leqslant C\norme{u}_{W^{m,p}}. \]
%\end{theorem}

\subsection{Parabolic estimates}

We recall estimates on the heat semigroup $\{e^{t\Delta}\}_{t\geq 0}$ in a bounded domain $\Omega$ with a homogeneous Neumann boundary condition. They can be found in \cite{CHKMR}, Lemma A.1 and \cite{HW}, Lemma 2.1.
\begin{lemma}\label{lem:LpLq}
	Let  $\lambda_1>0$ denote the first nonzero eigenvalue of $-\Delta$ in a bounded domain $\Omega\subset \R^d$ under the Neumann boundary conditions. For all  $1 \le q \le p \le +\infty$, there exist constants $C=C_1(p,q,\Omega)$ such that
	\begin{itemize}
		\item for all $f \in L^q (\Omega)$ satisfying $\int_\Omega f d x=0$ and all $t>0$,
		\begin{equation}\label{G1}
		\norme{{e^{t\Delta} f}}_p \le C\left(1+t^{- \frac{d}{2}\left( \frac{1}{q}-\frac{1}{p} \right)}\right)e^{-\lambda_1t}\norme{f}_q;
		\end{equation}
		\item  for all $f \in L^q (\Omega)$ and all $t>0$,
		\begin{equation}\label{G2}
		\norme{{e^{t\Delta} f}}_p \le C\big{(}1+ t^{- \frac{d}{2}\left( \frac{1}{q}-\frac{1}{p} \right)}\big{)}\norme{f}_q;
		\end{equation}
		
		\item  for all $f \in L^q (\Omega)$ and all $t>0$,
		\begin{equation}\label{G3}
		\norme{{\nabla \left(e^{t\Delta} f\right)}}_p \le C \left(1+t^{- \frac{d}{2}
			\left(\frac{1}{q}-\frac{1}{p}\right) - \frac{1}{2}}\right)e^{-{\lambda_1} t}\norme{f}_q;
		\end{equation}
		
		\item for all $f \in (L^q (\Omega))^d$ and all $t>0$,
		\begin{equation}\label{G4}
		\norme{{e^{t\Delta} \nabla \cdot f}}_p \le C \left(1+t^{- \frac{d}{2}
			\left(\frac{1}{q}-\frac{1}{p}\right) - \frac{1}{2}}\right)e^{-{\lambda_1} t}\norme{f}_q,
		\end{equation}
		provided $q>1$;
		\item  for all $f \in W^{1,p} (\Omega)$ and all $t>0$,
		\begin{equation}\label{G5}
		\norme{{\nabla e^{t\Delta} f}}_p \le C e^{-{\lambda_1} t}\norme{\nabla f}_p,
		\end{equation}
		provided $p\in [2,\infty)$.
	\end{itemize}
\end{lemma}

%
%\begin{theorem}\label{lem:LpLq}
%	Let $\lambda_1$ denote the first nonzero eigenvalue of $-\Delta$ in $\Omega\subset \R^d$ under the Neumann boundary conditions. For all $1\leqslant q \leqslant p \leqslant +\infty$, there exists a constant $C>0$ depending on $p,q,\Omega$ such that for all $f\in L^q(\Omega)$ and all $t>0$,
%	\[ \norme{\nabla \e^{t\Delta }f}_p \leqslant C t^{-\frac{d}{2}(\frac1q - \frac1p)-\frac12}\e^{-\frac{\lambda_1}{2}t}\norme{f}_q. \]
%\end{theorem}

\subsection{Other results}

One can find the following result and its proof in \cite{BF}, Theorem II.5.16.

\begin{lemma}[Aubin-Lions theorem]\label{Aubin-Lions}
	Let $X_0,X,X_1$ be three Banach spaces such that $X_0$ embeds compactly in $X$ and $X$ embeds continuously in $X_1$. Let $p,q\in[1,+\infty]$.
	We define
	\[  \mathcal{Y}=\{  u\in L^{p}(0,T;X_0) \ | \ \dfrac{\partial u}{\partial t}\in L^{q}(0,T;X_1) \},  \]
	where $\dfrac{\partial u}{\partial t}$ is the time derivative in the sense of distributions of $u$. Then,
	\begin{itemize}
		\item if $p<+\infty$, $\mathcal Y$ is compactly embedded in $L^{p}(0,T;X)$;
		\item if $p=+\infty$ and $q>1$ then $\mathcal Y$ is compactly embedded in $C([0,T],X)$;
	\end{itemize} 
\end{lemma}

One can find in \cite{M1}, Lemma 5.3 the following result:
\begin{lemma}[Gradient inequality in convex domains]\label{Conv}
	Let $\Omega\subset\R^d$ be a convex open set of smooth boundary. Let $u\in\mathscr{C}^2(\bar\Omega)$ such that $\dfrac{\partial u}{\partial \nu}=0$ on $\partial \Omega$. Then,
	\[ \dfrac{\partial}{\partial \nu} |\nabla u|^2 \leqslant 0 \qquad \mathrm{on} \ \partial\Omega. \]
\end{lemma}

We need in the proof of Lemma \ref{lemme2} the following result:
\begin{lemma}\label{upper_sol}
	Let $C>0$ and $w_0>0$. Consider the ordinary differential equation
	\begin{equation}\label{Some_ODE_bis} w'(t) =  C (1 + t + 2\sqrt{t}) \Big(1 + w(t)^\frac 3 p \left(1 + \log^+(w(t)) \right) \Big) w(t).   \end{equation}
	with initial condition $w(0)=w_0$. Then, there exists a time $T>0$ such that there exists a unique $C^1$ solution $w$ on $[0,T]$.
\end{lemma}

\noindent{\bf Proof.} 
	Since the equation can be recast in the form
	\[ w'(t) = f(t,w(t)),  \]
	with $f(t,w) =  C (1 + t + 2\sqrt{t}) \Big(1 + w^\frac 3 p \left(1 + \log^+(w)) \right) \Big) w$ a continuous function on $\R\times(0,+\infty)$ which is Lipschitz continuous with respect to $w$, the result is a consequence of the application of the Picard-Lindelöf theorem.
\qed
~\\

The following blow-up result for a particular ODE is widely known and used. However, we could not find any precise reference for it. This is why we present a proof below.

\begin{lemma}\label{BlowupODE}
	Let $\alpha_1,\alpha_2,\alpha_3 >0$ and $p\in(1,+\infty)$. Consider the ordinary differential equation
	\[  y'(t) = \alpha_1 y(t)^{1+\frac1p} - \alpha_2y(t) - \alpha_3  \]
	with initial condition $y(0)=y_0$. If $y_0$ is large enough with respect to $\alpha_1,\alpha_2,\alpha_3$ and $p$, then the unique local-in-time solution of this initial value problem blows-up in finite time at $T_{\mathrm{max}}$ and
	\[  \lim_{t\to T_{\mathrm{max}} } y(t) = +\infty.  \]
\end{lemma}

\noindent{\bf Proof.}
Let $\beta_0>0$. Denote for all $t>0$
\[ z(t) = y(t) - \beta_0. \]
If $z$ blows-up in finite time, then $y$ blows-up in finite time.
Assume $y_0$ is large enough so that $z(0)>0$ and for all $ x \geqslant y_0$,
\begin{equation}\label{cond_0}\alpha_1 x^{1+\frac 1 p} - \alpha_2x  - \alpha_3>0,\end{equation}
which is possible because this function tends to $+\infty$ when $x$ tends to $+\infty$. Then, for all $t>0$, $y'(t)>0$ and since $z'(t)=y'(t)$, it implies that $z'(t)>0$. Thus, $z$ is increasing.\\
We compute
\[
\begin{array}{rcl}
\displaystyle 	z'(t)&=&\alpha_1 y(t)^{1+\frac1p} - \alpha_2y(t) - \alpha_3\\
					 &=&\displaystyle  \alpha_1 (y(t) - \beta_0 + \beta_0)^{1+\frac1p} - \alpha_2(y(t) - \beta_0 + \beta_0) - \alpha_3\\
					 &\geqslant&\displaystyle \alpha_1 (y(t)-\beta_0)^{1+\frac1p} + \alpha_1(1+\frac1p)(y(t)-\beta_0)^\frac1p \beta_0 -\alpha_2 (y(t)-\beta_0) - (\alpha_2\beta_0 + \alpha_3)
\end{array}
\]
Hence, for all $t>0$,
\[ z'(t) \geqslant \alpha_1z(t)^{1+\frac1p} - \alpha_2z(t) +\alpha_1\beta_0(1+\frac1p)z(t)^\frac1p - \alpha_2\beta_0 - \alpha_3. \]
Since for all $t>0$, $z(t) > z(0)>0$, it follows that $\alpha_1\beta_0(1+\frac1p)z(t)^\frac1p  - \alpha_2\beta_0 - \alpha_3 > \alpha_1\beta_0(1+\frac1p)z(0)^\frac1p  - \alpha_2\beta_0 - \alpha_3$. If $z(0)$ is large enough we have
\begin{equation}\label{cond_1} \alpha_1\beta_0(1+\frac1p)z(0)^\frac1p  - \alpha_2\beta_0 - \alpha_3 > 0.\end{equation}
Then for all $t>0$, $z'(t) \geqslant \alpha_1z(t)^{1+\frac1p} - \alpha_2 z(t)$, and thus
\[  z'(t) \geqslant z(t) ( \alpha_1z(t)^\frac 1 p  - \alpha_2 ).  \]
Since $z(t)$ is increasing, if 
\begin{equation}\label{cond_2}\frac 1 2 \alpha_1z(0)^\frac 1 p - \alpha_2 > 0,\end{equation}
then we have that for all time $t>0$, $\frac 1 2 \alpha_1z(t)^\frac 1 p - \alpha_2 > 0$ and thus
\[  z'(t) \geqslant \frac 1 2 \alpha_1 z(t)^{1+\frac 1 p}.  \]
The equation $w'(t) = \frac 1 2 \alpha_1 w(t)^{1 + \frac 1 p}$ blows-up in finite time for all positive initial data.
Therefore, by the comparison principle for ordinary differential equations, if $y_0$ satisfies \eqref{cond_0} and is large enough to make $z(0)$ satisfy \eqref{cond_1} and \eqref{cond_2}, then $z$ blows-up in finite time; hence the result.
%  and then
% \[  y(t) \geqslant \dfrac{1}{\left(\frac{1}{ \varepsilon^\frac1p\beta_0^\frac1p }-\frac{t}{p}\right)^p} + (1+\varepsilon)\beta_0.\]
% The solution $y$ then blow-up in finite time and since it has the same derivative as $z$, it has the same blow-up rate.
\qed
~\\
%
%\begin{theorem}\label{lem:LpLq2}
%	For all $f\in W^{2,\infty}(\Omega)$,
%	\[ \norme{\nabla \e^{t\Delta} f}_\infty \leqslant \norme{f}_{W^{2,\infty}}.  \]
%\end{theorem}
%

\noindent{\Large\textbf{Acknowledgements}}\\
~\\ 
We would like to thank Laurent Moonens, Takayoshi Ogawa and Angela Stevens for useful discussions on some technical results.
%
%{\bf Proof.}
%
%\qed

%\newpage
%\bibliographystyle{siam}
%\bibliography{Biblio}

\begin{thebibliography}{99}



\bibitem{ADN}
	S. Agmond, A. Douglis, L. Nirenberg, \textit{Estimates near the boundary for solutions of elliptic partial differential equations satisfying general boundary conditions I},
	Communications on pure and applied mathematics, \textbf{XII}
	(1959), 623--727.


\bibitem{AMM}
A. Aotani, M. Mimura, T. Mollee, {\it A model aided understanding of spot pattern formation in chemotactic E. coli colonies}, Japan J. Indust. Appl. Math. {\bf 27} (2010), 5--22.


\bibitem{BBTW}
	N. Bellomo, A. Bellouquid, Y. Tao, M. Winkler, {\it Toward a mathematical theory of Keller–Segel models of pattern formation in biological tissues}, Mathematical Models and Methods in Applied Sciences, \textbf{25(09)} (2015), 1663-1763.



\bibitem{B1}
P. Biler, {\sc Singularities of Solutions to Chemotaxis Systems},  De Gruyter Series in Mathematics and Life Sciences, 6, De Gruyter,  2020.




\bibitem{BHLM1}
	A. Bonami, D. Hilhorst, E. Logak, M. Mimura, \textit{A free boundary problem arising in a chemotaxis model}, Free boundary problems, theory and applications, Pitman Res. Notes Math, \textbf{363} (1996), 368--373.


\bibitem{BHLM2}
	A. Bonami, D. Hilhorst, E. Logak, M. Mimura, \textit{Singular limit of a chemotaxis-growth model}, Adv. Differentials Equations, \textbf{6} (2001), 1173--1218.

\bibitem{BF}
F. Boyer, P. Fabrie, {\sc Mathematical Tools for the Study of the Incompressible Navier-Stokes Equations and
Related Models. Applied Mathematical Sciences}, New York: Springer, 2013.


\bibitem{B2}
H. Brezis, {\sc Functional Analysis, Sobolev Spaces and Partial Differential Equations}, Springer-Verlag New York, 2011.


\bibitem{BB1}
E. O. Budrene, H. Berg, {\it Complex patterns formed by motile cells of Escherichia coli}, Nature {\bf 349} (1991), 630--633.


\bibitem{BB2}
E. O. Budrene, H. Berg, {\it Dynamics of formation of symmetrical patterns by chemotactic bacteria}, Nature {\bf 376} (1995), 49--53.


\bibitem{C}
 X. Cao, {\it Global bounded solutions of the higher-dimensional Keller-Segel system under smallness conditions in optimal spaces}, Discrete \& Continuous Dynamical Systems-A, \textbf{35} (2015), 1891--1904.
%

\bibitem{C2}
X. Chen, \textit{Generation and propagation of interfaces in reaction-diffusion systems}, Transactions of the American Mathematical society, \textbf{334-2} (1992), 877--913.

\bibitem{CHKMR}
R. Celinski, D. Hilhorst, G. Karch, M. Mimura, P. Roux, \textit{Patterns in PDE model of chemotactic E. Coli colonies},
arXiv preprint arXiv:2003.06010

%
%\bibitem{CP}
%L. Corrias, B. Perthame, {\it Critical space for the parabolic-parabolic Keller–Segel model in Rd}.
%Comptes Rendus Mathematique  {\bf 342(10)} (2006), 745--750.
%

%
%
%\bibitem{CPZ}
%L. Corrias, B. Perthame, H. Zaag {\it A chemotaxis model motivated by angiogenesis}.
%Comptes Rendus Mathematique {\bf 336(2)} (2003), 141--146.



%
%
%\bibitem{D}
%	F. Dkhil, \textit{On the analysis of reaction-diffusion-advection systems arising in chemical and biomathematical models}, PhD-thesis, University of Cergy-Pontoise, France, (2002).



\bibitem{FHMW}
E. Feireisl, D. Hilhorst, M. Mimura, R. Weidenfeld, {\it On a nonlinear diffusion system with resource-consumer interaction}, Hiroshima Math. J. {\bf 33} (2003), 253--295.

\bibitem{FGM}
X. Fu, Q. Griette, P. Magal, \textit{A cell–cell repulsion model on a hyperbolic Keller–Segel equation}, Journal of Mathematical Biology, (2020), 1-44.

\bibitem{FMT}
	M. Funaki, M. Mimura, T. Tsujikawa, \textit{Travelling front solutions arising in a chemotaxis-growth model}, RIMS Kokyuroku, \textbf{1135} (2000), 52-76.


\bibitem{G}
D. Gilbarg, N. S. Trudinger, {\sc Elliptic Partial Differential Equations of Second Order}, Springer-Verlag, second edition, 2001.



\bibitem{HHS}
M. Henry, D. Hilhorst, R. Schätzle, \textit{Convergence to a viscosity solution for an advection-reaction-diffusion equation arising from a chemotaxis-growth model}, Hiroshima Math. J., \textbf{29} (1999), 591-630.

\bibitem{HP}
T. Hillen, K. J. Painter, \textit{A user’s guide to PDE models for chemotaxis}, Journal of mathematical biology, 58(1-2), (2009),  183.

\bibitem{H1}
D. Horstmann, \textit{From 1970 until present: the Keller-Segel model in chemotaxis and its consequences}, Jahresber. DMV {\bf 105} (2003), 103--165.


\bibitem{H2}
	D. Horstmann, \textit{From 1970 until present : the Keller-Segel model in chemotaxis and its consequences. II},
	 Max-Planck-Institut für Mathematik in den Naturwissenschaften Leipzig, \textbf{106}, (2004), 51--69


\bibitem{HW}
D. Horstmann, M. Winkler, {\it Boundedness vs. blow-up in a chemotaxis system}, J. Diff. Eqns., {\bf 215-1}, (2005), 52--107.



\bibitem{HMT}
P. P. Htoo, M. Mimura, I. Takagi, {\it Global solutions to a one-dimensional nonlinear parabolic system modeling colonial formation by chemotactic bacteria}, Adv. Stud. Pure Math. {\bf 47} (2007), 613--622.


%
%
%\bibitem{JL}
%W. J\"ager, S. Luckhaus, {\it On explosions of solutions to a system of partial differential equations modelling chemotaxis}, Trans. Amer. Math. Soc. {\bf 329-2} (1992), 819--824.




\bibitem{KS1}
	K. Kang, A. Stevens, \textit{Blow-up and global solutions in a chemotaxis-growth system}, Nonlinear Analysis, \textbf{135} (2016), 57--72.


\bibitem{KS2}
E. F. Keller, L. A. Segel, {\it Initiation of slime mold aggregation viewed as an instability}, J. Theoret. Biol. {\bf 26} (1970), 399--415.


%
%
%\bibitem{KS3}
%E. F. Keller, L. A. Segel, {\it Model for chemotaxis}, J. Theor. Biol. {\bf 30} (1971), 225--234.

%\bibitem{KS4}
%  H. Kozono, Y. Sugiyama, {\it Local existence and finite time blow-up of solutions in the 2-D Keller-Segel system}, J. evol. equ. {\bf 8} (2008), 353-378.


%\bibitem{KO}
%	M. Kurokiba, T. Ogawa, {\it Finite time blow-up of the solution for a nonlinear parabolic equation of drift-diffusion type}, Differential Integral Equations {\bf 16} (2003), 427--452.




\bibitem{L1}
	J. Lankeit, \textit{Eventual smoothness and asymptotics in a three-dimensional chemotaxis system with logistic source}, J. Differential Equations,
	\textbf{258} (2015), 1158--1191.


\bibitem{L2}
	J. Lankeit, \textit{Long-term behaviour in a chemotaxis-fluid system with logistic source}, Mathematical Models and Methods in Applied Sciences, \textbf{26} (2016), 2071--2109.


\bibitem{L3}
	J. Lankeit, \textit{Chemotaxis can prevent thresholds on population density},
	Discrete and continuous dynamical systems Series B,
	\textbf{20,5} (2015).


\bibitem{M1}
	H. Matano, \textit{Asymptotic behavior and stability of solutions of semilinear diffusion equations}, Publ. RIMS, Kyoto university,
	\textbf{15} (1979), 401--454.



%\bibitem{MSM}
%	M. Mimura, H. Sakaguchi, M. Matsushita, {\it Reaction-diffusion modelling of bacterial colony patterns}, Phys. A {\bf 282} (2000), 283--303.

\bibitem{MT}
M. Mimura, T. Tsujikawa, {\it Aggregating pattern dynamics in a chemotaxis model including growth}, Phys. A {\bf 230} (1996), 499--543.

%
%
%\bibitem{M2}
%	N. Mizoguchi,
%	{\it Finite-time blowup in Cauchy problem of parabolic-parabolic chemotaxis system.}
%	J. Math. Pures Appl. {\bf (9) 136} (2020), 203--238.

\bibitem{MY}
M. Mizukami, T. Yokota,  \textit{A unified method for boundedness in fully parabolic chemotaxis systems with signal‐dependent sensitivity}. Math. Nachrichten, \textbf{290} (2017), 2648--2660.


\bibitem{Private2}
L. Moonens, \textit{Private communication}.


%
%\bibitem{M3}
%X. Mora,
%{\it Semilinear Parabolic Problems Define Semiflows on $C^k$ Spaces.}
%Trans. Amer. Math. Soc. {\bf 278} (1983), 21--55.

%
%\bibitem{N1}
%T. Nagai, {\it Blowup of nonradial solutions to parabolic-elliptic systems modeling chemotaxis in two-dimensional domains}, J. Inequal. Appl. {\bf 6} (2001), 37--55.


%
%\bibitem{N2}
%G. Nardi, {\it Schauder estimate for solutions of Poisson’s
%	equation with Neumann boundary condition}, arXiv:1302.4103v2 (2013).


\bibitem{OT}
	T. Ogawa, Y. Taniuchi, \textit{On blow-up criteria of smooth solutions to the 3-D Euler equations in a bounded domain} Journal of Differential Equations, \textbf{190,1} (2003), 39--63.


\bibitem{Private}
T. Ogawa, \textit{Private communication}.


\bibitem{OTYM}
	K. Osaki, T. Tsujikawa, A. Yagi, M. Mimura, \textit{Exponential attractor for a chemotaxis-growth system of equations}, Nonlinear Anal. TMA, \textbf{51} (2002), 119--144.

%\bibitem{OY}
%	K. Osaki, A. Yagi, {\it Finite dimensional attractor for one-dimensional Keller-Segel equations}, Funkc. Ekvac. {\bf 44} (2001), 441--469.

\bibitem{PD}
B. Perthame,A. L. Dalibard, \textit{Existence of solutions of the hyperbolic Keller-Segel model}, Transactions of the american mathematical society, \textbf{361(5)} (2009), 2319-2335.


%
%\bibitem{PPLP}
%A. A. Polezhaev, R. A. Pashkov, A. I. Lobanov, I. B. Petrov, {\it Spatial patterns formed by chemotactic bacteria Escherichia coli}, Int. J. Dev. Biol. {\bf 50} (2006), 309--314.
%
%\bibitem{R}
%F. Rothe, {\sc Global solutions of reaction-diffusion systems}, Lecture Notes in Mathematics, vol. 1072. Springer-Verlag, Berlin-Heidelberg-New York-Tokyo, 1984.

%
%\bibitem{S1}
%	R. Serrano, \textit{An alternative proof of the Aubin-Lions lemma}, Arch. Math., \textbf{101} (2013), 253--257.
%
%\bibitem{SK}
%N. Shigesada, K. Kawasaki, {\it Modeling Complex Patterns in Bacterial Colonies}, Seibutsu Butsuri {\bf 40} (2000), 151--155 (in Japanese).



\bibitem{S}
E. Stein, T. S. Murphy, \textit{Harmonic analysis: real-variable methods, orthogonality, and oscillatory integrals. Vol. 3}, Princeton University Press (1993).


\bibitem{TW}
J. Tello, M. Winkler,
{\it A Chemotaxis System with Logistic Source},  Comm. in P.D.E., {\bf 32,6} (2007), 849--877.


%
%\bibitem{TLABCSR}
%L. Tsimring, H. Levine, I. Aranson, E. Ben-Jacob, I. Cohen, O. Shochet, W. N. Reynolds, {\it Aggregation patterns in stressed bacteria}, Phys. Rev. Lett. {\bf 75} (1995), 1859--1862.

%
%\bibitem{TLM}
%R. Tyson, S. R. Lubkin, J. D. Murray, {\it Model and analysis of chemotactic bacterial patterns in a liquid medium}	, J. Math. Biol. {\bf 38} (1999), 359--375.


\bibitem{W1}
M. Winkler,	{\it How Far Can Chemotactic Cross-diffusion Enforce Exceeding Carrying Capacities?}, J. Nonlinear Science, {\bf 24}	(2014), 809--855.


\bibitem{W2}
M. Winkler, \textit{Aggregation vs. global diffusive behavior in the higher-dimensional Keller–Segel model}, J.  Differential Equations, \textbf{248(12)} (2010), 2889--2905.

\bibitem{W3}
M. Winkler, \textit{Blow-up in a higher-dimensional chemotaxis system despite logistic growth restriction}, J. Math. Analysis Appl., \textbf{384} (2011), 261--272.
%
%\bibitem{W4}
%	M. Winkler, {\it 
%		Finite-time blow-up in the higher-dimensional parabolic-parabolic Keller-Segel
%		system.}
%	J. Math. Pures Appl. {\bf (9) 100} (2013),  748--767.


%\bibitem{W5}
%M. Winkler, {\it Absence of collapse in a parabolic chemotaxis system with signal-dependent sensitivity}, Math. Nachr. {\bf 283} (2010), 1664--1673.


%\bibitem{W6}
%M.\ Winkler,\ {\it Aggregation vs. global diffusive behavior in the higher-dimensional Keller-Segel model},\ J.\ Differential\ Equations {\bf 248} (2010), 2889--2905.


%
%\bibitem{W7}
%I. Wood,	{\it Maximal $L^p-$regularity for the Laplacian on Lipschitz domains}, Mathematische Zeitschrift {\bf 255}	(2007),  855–875.
%
%
%\bibitem{X}
% T. Xiang, \textit{Global dynamics for a diffusive predator–prey model with prey-taxis and classical Lotka–Volterra kinetics}. Nonlinear Analysis: Real World Applications, \textbf{39} (2018), 278--299.

%
%\bibitem{Y}
%A. Yagi, {\sc Abstract parabolic evolution equations and their applications}, Springer Monographs in Mathematics. Springer-Verlag, Berlin, 2010




\end{thebibliography}

\end{document}